\documentclass[12pt]{article}

\def\qed{\quad \vrule height7.5pt width4.17pt depth0pt} %Ed's QED

\newtheorem{theorem}{Theorem}
\newtheorem{lemma}{Lemma}
\newtheorem{cor}{Corollary}
\newtheorem{prop}{Proposition}

\begin{document}

\title{Hitting properties of a random string}
\author{Carl Mueller$^{1}$ and Roger Tribe\\\\
\\Dept. of Mathematics  \\University of Rochester \\Rochester, NY  14627, USA
\\E-mail:  cmlr@troi.cc.rochester.edu
\\
\\Mathematics Institute \\University of Warwick \\Coventry  CV4 7AL, UK
\\E-mail: tribe@maths.warwick.ac.uk}
\date{}
\maketitle

\footnotetext[1]{Supported by an NSF travel grant and an NSA grant. 


AMS 1991 {\em subject classifications} Primary, 60H15; Secondary, 35R60, 
35L05.}

%\newpage

\begin{abstract}
We consider Funaki's model of a random string taking values in 
$\mathbf{R}^d$.  It is specified by the following stochastic PDE,
\[
\frac{\partial u(x)}{\partial t}=\frac{\partial^2 u(x)}{\partial x^2} 
+\dot{W}.
\]
where $\dot{W}=\dot{W}(x,t)$ is two-parameter white noise, also taking 
values in $\mathbf{R}^d$.  We find the dimensions in which the string hits 
points, and in which it has double points of various types.  We also study 
the question of recurrence and transience.
\end{abstract}
%
%\newpage
%
%%%%%%%%%%%%%%%%%%%%%%%%%%%%%%%%%%%%%%%%%%%%%%%%%%%%
\section{Introduction} \label{s1}
\setcounter{equation}{0}
%
%%%%%%%%%%%%%%%%%%%%%%%%%%%%%%%%%%%%%%%%%%%%%%%%%%%%
In this paper we study hitting problems, double points  and recurrence 
questions for the following model
of a random string, first introduced by Funaki \cite{fun84}:
\begin{equation}
\label{spde}
\frac{\partial u(x)}{\partial t}=\frac{\partial^2 u(x)}{\partial x^2} 
+\dot{W}.
\end{equation}
Here $ \dot{W}(x,t) $  is a $\mathbf{R}^d$ valued space-time white noise and 
$(u_t(x):t \geq 0, \, x \in \mathbf{R})$ is a continuous $\mathbf{R}^{d}$ 
valued process.  We give details of the meaning of this equation below.  We 
also consider the analogous random loop, that is a solution  indexed over 
$x \in \mathbf{T} = \mathbf{R}$ (mod 1),  the circle.  We give a brief 
motivation for this model.  Under Newton's law of motion, the equation for 
the damped motion of a particle of mass $m$ in a force field $F$ is 
\[
m \frac{\partial^2 x(t)}{\partial t^2}
        =F(x(t)) - \lambda \frac{\partial x(t)}{\partial t}.
\]
However if the particle has low mass, or the force field and damping are 
strong, then the motion is well approximated by Aristotle's law
\[
\frac{\partial x(t)}{\partial t}= \lambda^{-1} F(x(t)),  
\]
which says that the velocity is proportional to the force.  In the same way,
the usual equation for an elastic string can be approximated by a heat 
equation.  Allowing the string to move in $\mathbf{R}^d$ leads us to look
for $\mathbf{R}^d$-valued solutions.  If  the string is influenced by white 
noise we arrive at the  equation (\ref{spde}).  Simple linear scaling allows 
us to set all parameters to one. 

Before proceeding further, we now give an outline of our main results, 
which are stated in greater detail in theorems 1,2 and 3, later in the 
paper.  We say that the random string  hits a point $z \in \mathbf{R}^d$ if 
$u_t(x) = z$ for some $t>0, \, x \in \mathbf{R}$.  We shall show the 
following properties hold, each as an  almost sure event:
\begin{itemize}
\item The random string hits points if and only if $d<6$;
\item For fixed $t_0>0$, there exist points $x,y$ such that 
$u_{t_0}(x)=u_{t_0}(y)$ if and only if $d<4$.  
\item There exist points $(t,x)$ and $(t,y)$ such that $u_t(x)=u_t(y)$ if 
and only if $d<8$;
\item There exist points $(t,x)$ and $(s,y)$ such that $u_t(x)=u_s(y)$ if 
and only if $d<12$.
\end{itemize}
The string is H\"{o}lder continuous of any order less than $1/2$ in space 
and $1/4 $ in time.  This suggests that that the range of the process 
$(t,x) \to u_t(x)$ might be 6 dimensional
and leads to the guess that dimension $d=6$ is the critical for hitting 
points and $d=12$ is critical for double points.  The second assertion in 
the above list is obvious, once we show that a for certain version of the 
string, $x\to u_{t_0}(x)$ is a Brownian motion parameterized by $x$.  (See 
Remark 1 after Theorem \ref{multiple-points}, at the beginning of Section 
\ref{s5}).  It is well-known that Brownian motion has double points if and 
only if $d<4$.

Thanks to the Markov property and potential theory, we have detailed
information about the hitting behavior of Brownian motion and other stochastic
processes.  Much less is known about the hitting behavior of random fields. 
There are 2 prominent exceptions to this statement.  The first involves 
random fields of the form
\[
Z(t_1,\ldots,t_n)=(X^{(1)}_{t_1},\ldots,X^{(n)}_{t_n})
\]
where $X^{(1)}_{t_1},\ldots,X^{(n)}_{t_n}$ are independent processes.  See
Fitzsimmons and Salisbury \cite{fs89} for this kind of work.  The second case
involves the Brownian sheet and other random fields with the multi-parameter
Markov property (see Orey and Pruitt \cite{op73}, Hirsch and Song \cite{hs95},
and Khoshnevesen and Shi \cite{ks99}).  In both of these cases, there is good
information about what kind of sets the process will hit with positive
probability.  The above examples rely heavily on the Markov properties for the
random fields in question.  Peres \cite{per96} has also done some beautiful
work with applications to the hitting properties of random fields.  In
addition to hitting questions, Orey and Pruitt \cite{op73} also studied the
recurrence and transience of the Brownian sheet.

We employ two main methods in the proof.  The finiteness of hitting times 
and the existence of double points of various types is deduced, below the 
critical dimensions, by a straightforward inclusion-exclusion argument.  The 
more interesting direction is the proof of the non-existence of such events 
in critical dimensions.  For this we make use of a stationary pinned version 
of the random string, constructed in section
\ref{s2}.  Starting the string off as a two sided Brownian motion leads to a 
solution for which the distribution of 
$(u_{t+t_0}(x)-u_{t_0}(0): t \geq 0, \, x \in \mathbf{R})$ does not depend 
on $t_0$.  The word `pinned' refers to looking at 
the image of the string under the map $f \to f - f(0)$.  The image under 
this map is still a Markov process 
and the law of a two sided Brownian motion is its unique stationary 
distribution.  The stationary pinned string has a simple scaling property 
that allows us to use scaling arguments, for example to show the Lebesgue 
measure of the range is zero in suitable critical dimensions.  An absolute 
continuity argument, given in section \ref{s3},  will show that our results 
hold for the string, not just in its stationary version, but also for the 
random loop.  Sections \ref{s4} and \ref{s5} contain the arguments for hitting 
points and for double points respectively.  It makes sense to ask for 
recurrence properties of the stationary string,  in the same spirit as Orey  
and Pruitt's results for the Brownian sheet.  In section \ref{s6} we show 
that the random string is recurrent if and only if $d \leq 6$.  We finish 
this introduction by briefly discussing existence of solutions to 
(\ref{spde}), giving a simple inclusion-exclusion type lemma, and 
introducing common notations used in the text. 

The components $ \dot{W}_1(x,t),\ldots, \dot{W}_d(x,t)$ of the vector noise 
$\dot{W}(x,t)$ are 
independent space-time white noises, which are generalized Gaussian processes
with covariance given by
$ E\left[\dot W_i(x,t)\dot W_i(y,s)\right]=\delta(t-s)\delta(x-y) $.  That 
is, $W_i(f)$ is a random field indexed by functions 
$f\in\mathbf{L}^2([0,\infty)\times\mathbf{R})$, and for two such test 
functions $f,g\in\mathbf{L}^2([0,\infty)\times\mathbf{R})$
we have 
\[
E\left[W_i(f)W_i(g)\right]=\int_{0}^{\infty}\int f(t,x)g(t,x)dxdt.
\]
Heuristically,
\[
W_i(f) = \int_{0}^{\infty}\int f(t,x)W(dx\,dt)
\]
We suppose that the noise is adapted with respect to a filtered probability 
space $(\Omega, {\cal F}, ({\cal F}_t), P)$, where ${\cal F}$ is complete and 
$({\cal F}_t)$ is right continuous, in that $W(f)$ is 
${\cal F}_t$-measurable whenever $f$ is supported in $[0,t] \times \mathbf{R}$.  

The initial conditions play an unimportant role in the properties we study 
for the solutions to
(\ref{spde}).  We may take any initial conditions that are suitable for the 
deterministic heat equation.  To be concrete, and so that we may apply 
results from the literature, we shall take initial conditions in 
$\mathcal{E}_{exp}$, the space of  continuous functions of at most 
exponential growth,
defined by $\mathcal{E}_{exp} = \cup_{\lambda >0} \mathcal{E}_{\lambda}$ where
\[
\mathcal{E}_{\lambda} = \{ f \in C(\mathbf{R}, \mathbf{R^d}):  
\mbox{ $|f(x)| \exp(-\lambda |x|) \to 0 $ as $x \to \pm \infty$} \}.
\]
We define a solution to (\ref{spde}) to be a $({\cal F}_t)$ adapted,  
continuous random field
 $(u_t(x): t \geq 0, \, x \in \mathbf{R})$ satisfying 
\begin{itemize}
\item[(i)] $u_0 \in \mathcal{E}_{exp}$ almost surely and is adapted to 
${\cal F}_0$, 
\item[(ii)] for each $t>0$ there exists $\lambda>0$ so that 
$u_s \in \mathcal{E}_{\lambda}$, 
for all $s \leq t$ almost surely, 
\item[(iii)] for each $t>0$ and $x \in \mathbf{R}$, the following Green's 
function representation holds
\begin{equation} \label{solution}
u_t(x)=\int G_t(x-y)u_0(y)dy + \int_{0}^{t}\int G_{t-r}(x-y)W(dy\,dr).
\end{equation}
\end{itemize}
Here $G_t(x)$ is the fundamental solution  
$ G_t(x)=(4\pi t)^{-1/2}\exp (-x^2/4t) $.  We note that for initial 
conditions $u_0$ that are deterministic, or that are Gaussian fields 
independent of ${\cal F}_0$, the solutions are Gaussian fields.  For any 
deterministic initial condition in $\mathcal{E}_{exp}$ there is a version 
of the solution (\ref{solution}) satisfying the regularity condition (ii), 
and the laws of these solutions form a Markov family in time.

For $\phi: \mathbf{R} \to \mathbf{R}^d$, we write $(u_t,\phi)$ for the 
integral $\int u_t(x) \phi(x) dx$, whenever this is well defined.  The above 
definition of solutions is equivalent to a `weak' formulation, in that 
condition (iii) may be replaced by the following: for all $\phi$ smooth and 
of compact support
\[
 (u_t, \phi) = (u_0, \phi) + \int^t_0 (u_s, \Delta \phi) ds 
         + \int^t_0 \int \phi(y) W(dy\,ds).
\]
Pardoux \cite{par93} or Walsh \cite{wal86} are references for the basic 
properties of SPDEs driven by space-time white noise as used above.  The 
equivalence of the weak formulation is shown, in the case of real valued 
stochastic PDE solutions, in Shiga (\cite{shi94}).

We shall make frequent use of the following inclusion-exclusion type lemma. 
\begin{lemma}
\label{lower-bound}
Suppose that $ (A_i:i=1,\ldots,n)$ are events and $A=\bigcup_{i=1}^n A_i$ Then
\begin{equation} 
P(A) \geq 
 \frac{\big[\sum_{i=1}^{n}P(A_i)\big]^2}
  {\sum_{i=1}^{n}P(A_i)+2\sum_{1\le i<j\le n}P(A_i\cap A_j)}.  
\end{equation}
\end{lemma}
This lemma is similar in spirit to the standard inclusion-exclusion bound
\[ 
P(A) \geq \sum_{i=1}^{n}P(A_i)-\sum_{1\leq i < j \leq n}P(A_i\cap A_j).  
\]
For both lower bounds, one must find a lower bound for 
$\sum_{i=1}^{n}P(A_i)$.  But to obtain a useful lower bound using Lemma 
\ref{lower-bound}, one often only need show that 
$\sum_{1 \leq i < j \leq n}P(A_i \cap A_j)$ is comparable to 
$(\sum_{i=1}^{n}P(A_i))^2$, while for the standard inclusion-exclusion 
bound one must show that it is strictly less.  The well  known proof is 
easy, but so short that we include it.  For a real variable  $Z \geq 0$ we 
have 
\[
E[Z] = E[Z\mathbf{1}(Z>0)] \leq \left[E(Z^2) \, P(Z>0) \right]^{1/2}
\]
by the Cauchy-Schwartz inequality.  The lemma follows from this inequality, 
after rearranging, by taking $Z = \sum_i \mathbf{1}(A_i)$ and noting 
$P(A) = P(Z>0)$. 

Finally, here is some notation.  For $x \in \mathbf{R}^d$ and $r \geq 0$ we 
write $B_r(x)$ for  the box 
$\{y \in \mathbf{R}^d: |y_i-x_i| < r\}$.  For use in our inclusion-exclusion 
arguments we define space and time  grids of points by 
\begin{equation}
\label{grid}
\mbox{$t_{i,n} = i2^{-4n}$ and 
$ \quad x_{j,n} = j2^{-2n}$ for $n,\,i,\, j \in \mathbf{Z}$}.  
\end{equation}
Throughout the paper $C(d, T, \ldots)$ will denote a constant, whose 
dependency will be listed, but whose value is unimportant and may change 
from line to line.  Constants $c_1,c_2, \ldots $ with  subscripts denote 
specific constants which do not change and may be referred to later in the 
text.
%
%%%%%%%%%%%%%%%%%%%%%%%%%%%%%%%%%%%%%%%%%%%%%%%%%%%%
\section{The stationary pinned string} \label{s2}
\setcounter{equation}{0}
%
%%%%%%%%%%%%%%%%%%%%%%%%%%%%%%%%%%%%%%%%%%%%%%%%%%%
%
The motivation for the pinned string comes from the following calculation. 
Starting from zero initial conditions, the solution to (\ref{spde}) is given 
by $u_t(x) =  \int_0^t  \int G_{t-r}(x-z)W(dz\,dr)$.
The variance of the first component is given by 
\begin{eqnarray*}
E\left[\left( u^{1}_t(x) \right)^2 \right] 
&=& \int_{0}^{t} \int G_{r}(x-z)^2 dzdr \\ 
&=& \frac{1}{\sqrt{8 \pi}}  \int_{0}^{t} r^{-1/2}dr  \\
&=& \frac{1}{\sqrt{2\pi}}t^{1/2}
\end{eqnarray*}
and diverges to infinity as $t \to \infty$.  However the variance of a 
spatial increment has the following limit as $t\to\infty$.  
\begin{eqnarray*}
E\left[\left( u^{1}_t(x) - u^{1}_t(y)  \right)^2 \right] 
&=& \int_{0}^{t} \int (G_{r}(x-z) - G_r(y-z)) ^2 dzdr  \\
& \to & \int_{0}^{\infty} \int (G_{r}(x-z) - G_r(y-z)) ^2 dzdr \\
& = & |x-y|.
\end{eqnarray*}
One way to calculate the double integral above and justify the final 
equality is to apply Plancherel's theorem, using the Fourier transform 
$ \hat{G}_r(\theta) = (2 \pi)^{-1/2} \exp (-r \theta^2)$ of $G_r$, to 
rewrite this double integral as
\begin{eqnarray*}
 \lefteqn{ \frac{1}{2 \pi} \int^{\infty}_0 \int \exp(-2r \theta^2) 
 |\exp(i(x-y)\theta)-1|^2 d \theta dr} \\
&=&  \frac{|x-y|}{2 \pi} \int^{\infty}_0 \int \exp(-2s \eta^2) |\exp(i 
\eta)-1|^2 d \eta ds \\
& = & \frac{|x-y|}{4 \pi} \int \frac{|\exp(i \eta)-1|^2 }{\eta^2} d \eta,
\end{eqnarray*}
using the substitutions $\eta = \theta (x-y)$ and $s=r |x-y|^{-2}$.  This 
shows the answer is of the form 
$c_0 |x-y|$ and the remaining integral can be evaluated, for example by 
contour integration, showing that $c_0=1$.  The limiting variance $|x-y|$ is 
exactly that of a  two sided Brownian motion.  The idea is to start a 
solution with this covariance structure and check that the spatial 
increments are stationary in time.

Motivated by the calculation above, we take an initial function
$(U_0(x):x \in \mathbf{R})$ which is a two-sided $\mathbf{R}^d$-valued 
Brownian motion  satisfying $U_0 (0)=0$ and 
$E[ (U_0(x) - U_0(y) )^2] = |x-y|$, and which is independent of the white 
noise $\dot{W}$.  This can be created in the following way: take an 
independent space-time white noise $\tilde{W}$ and let 
\[
U_0(x) = \int^{\infty}_0 \int \left( G_r(x-z) - G_r(z) \right) 
\tilde{W}(dz\,dr).
\]
Indeed, the calculation above shows that this integral has the correct 
covariance.  We may assume, by extending the probability space if needed, 
that $U_0$ is ${\cal F}_0$-measurable.  The solution to  (\ref{spde}) driven 
by the noise $W(x,s)$ is then given by
\begin{eqnarray}
\label{stationarystring}
U_t(x) & = & \int  G_{t}(x-z)  U_0(z)  dz + \int^t_0  \int  G_r (x-z) W(dz\,dr) \\
& = & \int^{\infty}_0 \int \left( G_{t+r}(x-z) - G_{t+r}(z) \right) 
\tilde{W}(dz\,dr)  \nonumber\\
&& \hspace{.4in}  + \int^t_0  \int  G_r (x-z) W (dz\,dr). \nonumber
\end{eqnarray}
We call a continuous version of this process the stationary pinned string. 
Note that the components 
$(U^{i}_t(x):t \geq 0, \, x \in \mathbf{R})$ for $i=1,\ldots,d$ are 
independent and identically distributed. 
\begin{prop}
The components $(U^{(i)}_t(x): x \in \mathbf{R}, t \geq 0)$ of the 
stationary pinned string are mean zero Gaussian fields with the following 
covariance structure: for $x,y \in \mathbf{R}$, $t \geq 0$ 
\begin{equation} \label{spatialcovariance}
E\left[ \left(U_t^{(i)}(x)-U_t^{(i)}(y)\right)^2\right]  = |x-y|,
\end{equation}
and for  $x,y \in \mathbf{R}$, $0 \leq s < t$  
\begin{equation} \label{generalcovariance}
E \left[ \left(U^{(i)}_t(x)-U^{(i)}_s(y)\right)^2 \right] = 
(t-s)^{1/2}F\left( |x-y| (t-s)^{-1/2} \right)
\end{equation}
where
\[
F(a)=(2\pi)^{-1/2}+\frac{1}{2} \int\int G_1(a-z) G_1(a-z') \left( |z| + 
|z'| - |z-z'| \right) dz dz'.  
\]
$F(x)$ is smooth function, bounded below by $(2 \pi)^{-1/2}$, 
%satisfying $F(0)=(\sqrt{8}-1)/\sqrt{2 \pi }$ 
and $F(x)/|x| \to 1 $ as $|x| \to \infty$.  Furthermore there exists $c_1>0$ 
so that for all $x,y \in \mathbf{R}$, $0 \leq s \leq  t$ 
\begin{equation}
\label{boundoncovariance}
c_1 \left( |x-y| + |t-s|^{1/2} \right) \leq E \left[ 
\left(U^{(i)}_t(x)-U^{(i)}_s(y)\right)^2 \right]
\leq  2 \left( |x-y| + |t-s|^{1/2} \right).
\end{equation}
\end{prop}
\textbf{Proof.} 
Aiming for (\ref{spatialcovariance}), a simple calculation using the 
isometry for stochastic integrals, and the independence of the two 
integrals in (\ref{stationarystring}), gives
\begin{eqnarray*}
\lefteqn{ E\left[ \left(U_t^{(i)}(x)-U_t^{(i)}(y)\right)^2\right]   }\\
& = & E\left[ \left( \int^{\infty}_0 \int \left( G_{t+r}(x-z) - 
G_{t+r}(y-z) \right) \tilde{W}(dz\,dr) 
\right)^2 \right] \\
& & \hspace{.4in} + E\left[ \left( \int^t_0  \int  \left( G_r (x-z) - 
G_r(y-z) \right) 
 W (dz\,dr) \right)^2 \right] \\
& = & \int^{\infty}_0 \int \left( G_r (x-z)- G_r (y-z) \right)^2 dz dr\\
& = & |x-y|.
\end{eqnarray*}
To calculate (\ref{generalcovariance}) we use the fact that
\[
 U_t(x) = \int G_{t-s}(x-z) U_s(z) dz + \int^t_s \int G_{t-r}(x-z) W(dz\,dr)
\]
so that
\begin{eqnarray*}
\lefteqn{ E \left[ \left. \left(U^{(i)}_t(x)-U^{(i)}_s(y) \right)^2 \right|  
{\cal F}_s \right]} \\
& = & E \Bigg[ \bigg( \int G_{t-s}(x-z) (U^{(i)}_s(z)-U^{(i)}_s(y))dz \\
&&\hspace{.4in}
+ \int^t_s \int G_{t-r}(x-z) W(dz\,dr) \bigg)^2 \bigg| {\cal F}_s \Bigg] \\
& = & \left( \int G_{t-s}(x-z) 
(U^{(i)}_s(z)-U^{(i)}_s(y))dz \right)^2 + \int^t_s \int G^2_{t-r} (x-z) dz 
dr \\
& = & \left( \int G_{t-s}(x-z) (U^{(i)}_s(z)-U^{(i)}_s(y))dz \right)^2
+ \left( \frac{|t-s|}{2 \pi} \right)^{1/2}.
\end{eqnarray*}
Using (\ref{spatialcovariance}) we have
\begin{eqnarray*}
\lefteqn{  E \left[ \left( \int G_{t-s}(x-z)  \Big(U^{(i)}_s(z)-U^{(i)}_s(y)\Big) 
dz \right)^2 \right] } \\
& = & \int G_{t-s}(x-z) \int G_{t-s}(x-z')   \\
&&\hspace{.4in}\cdot E \bigg[\Big(U^{(i)}_s(z)-U^{(i)}_s(y)\Big)
 \Big(U^{(i)}_s(z')-U^{(i)}_s(y)\Big) \bigg] dz dz' \\
& = & \frac12 \int G_{t-s}(x-z) \int G_{t-s}(x-z') \\
&& \hspace{.4in}\cdot E \bigg[\Big(U^{(i)}_s(z)-U^{(i)}_s(y)\Big)^2  +  
\Big(U^{(i)}_s(z')-U^{(i)}_s(y)\Big)^2   \\
&&\hspace{.7in}- \Big(U^{(i)}_s(z)-U^{(i)}_s(z')\Big)^2   \bigg] dz dz' \\
& = & \frac12 \int G_{t-s}(x-z) \int G_{t-s}(x-z') \left( |z-y| + |z'-y| - 
|z-z'| \right) dz dz'.
\end{eqnarray*}
The scaling $ G_r(x) = r^{-1/2} G_1(x/r^{1/2}) $ now leads to the 
covariance formula (\ref{generalcovariance}).  The function $F(a)$ can be 
expressed in terms of exponentials and Gaussian error functions.  However,
the form given makes it clear that $F(a) \geq (2 \pi)^{-1/2}$.  The 
calculations needed to establish  the other properties of  $F(a)$ are 
straightforward and omitted.

The upper bound in (\ref{boundoncovariance}) follows directly from 
(\ref{spatialcovariance}) and (\ref{generalcovariance}).  This upper bound  
implies that there is a continuous 
version of $(t,x) \to U_t(x)$ and that this version satisfies the growth 
estimates (ii) in the definition of a solution (\ref{solution}).   The 
lower bound in (\ref{boundoncovariance}) is immediate from 
(\ref{spatialcovariance}) in the case $s=t$.  If $ |x-y| < |t-s|^{1/2}$ then 
we argue that
\begin{eqnarray*}
 E \left[ \left(U^{(i)}_t(x)-U^{(i)}_s(y)\right)^2 \right] & = & 
 |t-s|^{1/2} F(|x-y| \, |t-s|^{-1/2})  \\
& \geq & \frac{1}{\sqrt{2 \pi}}  |t-s|^{1/2}\\
& \geq & \frac{1}{2 \sqrt{2 \pi}}   \left( |x-y| + |t-s|^{1/2} \right) . 
\end{eqnarray*}
If $ 0< |t-s|^{1/2} \leq |x-y| $ then we argue that 
\begin{eqnarray*}
 E \left[ \left(U^{(i)}_t(x)-U^{(i)}_s(y)\right)^2 \right] &=&  |t-s|^{1/2} 
 F(|x-y| \, |t-s|^{-1/2}) \\
&  \geq  & |x-y|  \inf \{ F(z)/z: z \geq 1 \} \\
& \geq & \frac12  \inf \{ F(z)/z: z \geq 1 \}  \left( |x-y| + |t-s|^{1/2} 
\right) . 
\end{eqnarray*}
Combining the three cases, we may take 
\[
c_1 = \min \left(  (8 \pi)^{-1/2}, \inf \{F(z)/2z: z \geq 1\}  \right).   
\]
This finishes the proof.  \qed
\begin{cor}
The stationary pinned string has the following properties:
\begin{itemize}
\item {\bf Translation invariance} For any $t_0 \geq 0$ and 
$x_0 \in \mathbf{R}$ the field
\[
(U_{t_0 +t}(x_0 \pm x)- U_{t_0}(x_0): x \in \mathbf{R}, t \geq 0)
\]
 has the same law as the stationary pinned string.
\item {\bf Scaling} For $L>0$ the field 
\[
 ( L^{-1} U_{L^4 t}(L^2 x) : x \in \mathbf{R}, t \geq 0) 
\]
has the same law as the stationary pinned string.
\item {\bf Time reversal} For any $T>0$ the field
\[ 
(U_{T-t}(x)-U_T(0): x \in \mathbf{R}, \, 0 \leq t \leq T)
\]
has the same law as the stationary pinned string over the interval $[0,T]$.
\end{itemize} 
\end{cor}
\textbf{Proof.} The covariance formulae (\ref{spatialcovariance}) and 
(\ref{generalcovariance}),
together with $U_0(0)=0$, characterize the law of the stationary pinned 
string.  The translation invariance, scaling and time reversal follow 
immediately by checking  that this covariance structure is preserved.   \qed
%
%%%%%%%%%%%%%%%%%%%%%%%%%%%%%%%%%%%%%%%%%%%%%%%%%%%%
\section{Absolute Continuity Results} \label{s3}
\setcounter{equation}{0}
%
%%%%%%%%%%%%%%%%%%%%%%%%%%%%%%%%%%%%%%%%%%%%%%%%%%%%
There are general criteria for the absolute continuity of Gaussian random 
fields
(see for example Ibragimov and Rozanov \cite{ir78}).  For our differential 
equation setting we 
found it easier to exploit Girsanov's theorem (see Dawson \cite{daw78}
and Nualart and Pardoux \cite{np94} for applications to stochastic PDEs). 
The following lemma deals with solutions to the perturbed equation
\begin{equation} \label{perturbation}
\frac{\partial v_t(x)}{\partial x}= \frac{\partial^2 v_t(x)}{\partial x^2} 
+  h_t(x) + \dot W(x,t).
\end{equation}
where $h_t(x): [0,\infty) \times \mathbf{R} \to \mathbf{R}^d$  is an 
adapted, continuous function.  Solutions to (\ref{perturbation}) are defined 
as in the introduction, with an extra drift term; for example the weak 
formulation has the extra integral $\int^t_0 (h_s,\phi) ds$. 
\begin{lemma} \label{girsanov}
Suppose $(u_t(x))$ is a solution to (\ref{spde}) and $(v_t(x))$ is a 
solution to (\ref{perturbation}).  Suppose also that they have the same 
deterministic initial condition $u_0(x)=v_0(x) = f \in \mathcal{E}_{exp}$. 
Then either of the following two conditions on $h_t(x)$ is sufficient to 
imply that the laws $P_u^{(T)}$ and $P_v^{(T)}$ of the solutions $(u_t(x))$ 
and $(v_t(x))$, on the region $ (t,x) \in [0,T] \times \mathbf{R}$, are 
mutually absolutely continuous. 
\begin{itemize}
\item[(a)] The drift $h_t(x)$ is deterministic and satisfies 
$ \int_{0}^{T} \int |h_t(x)|^2 dx dt <\infty $.
\item[(b)] The drift $h_t(x)$ has compact support $A$ and is independent of 
$(W(dx\,dt): (t,x) \in A)$.
\end{itemize} 
\end{lemma}
The lemma is an  easy consequence of  Girsanov's change of measure theorem. 
Indeed suppose $(v_t(x))$ is a solution to (\ref{perturbation}), started at 
$f$, on the probability space $(\Omega, {\cal F}, ({\cal F}_t), P)$.  Define
\begin{equation} \label{Rnderivative}
\frac{dQ}{dP} = \exp \left( \int^T_0 \! \int h_t(x) \cdot W(dx\,dt) - 
\frac{1}{2} \int^T_0 \! \int
|h_t(x)|^2 dx \, dt \right).
\end{equation}
In the case $h_t(x)$ is deterministic the stochastic integral is Gaussian 
and the exponential defines a martingale.  Then under $Q$ the process 
$v_t(x)$ is a solution to (\ref{spde}) started at $f$ with respect to a new 
noise defined by $\tilde{W}(f) = W(f) - \int^t_0 (h_s,f) ds$.  This is 
easiest to check using the weak formulation of the equation and Levy's 
characterization of a space-time white noise $W$ (see Walsh \cite{wal86} 
Chapter 3).  In a similar way one can obtain a solution of 
(\ref{perturbation}) starting from a solution to (\ref{spde}).  

In case (b), the same proof works once one knows that the exponential in 
(\ref{Rnderivative}) defines
a true martingale.  Since $h_t(x)$ is continuous and adapted the stochastic 
integral in (\ref{Rnderivative}) is well defined and the formula for $dQ/dP$ 
defines a positive supermartingale.  It is sufficient to check that it has 
expectation $1$ to ensure it is a martingale.  Let $\mathcal{G}$ be the 
$\sigma$-field generated by
$( W(dx\,dt): (t,x) \in A)$.  Conditioned on $\mathcal{G}$, the stochastic 
integral is Gaussian, and so
\[
E  \left[ \left.\exp\left(\int_0^{T}\int h_t(x)\cdot W(dx\,dt) 
-\frac{1}{2}\int_0^{T}\int |h_t(x)|^2dxdt\right)
\right| \mathcal{G} \right] =1.
\]
Taking a further expectation shows the exponential has expectation $1$.

One consequence of the absolute continuity is that, under the conditions of 
the lemma,  solutions to (\ref{perturbation}) are unique in law and satisfy 
the Markov property. 
\begin{cor} \label{absolutecorollary1}
Suppose $(u_t(x))$ and $(\tilde{u}_t(x))$ are both solutions to 
(\ref{spde}).  For compact sets
$A \subseteq (0,\infty) \times \mathbf{R}$ the laws of the fields 
$(u_t(x): (t,x) \in A)$ and
$(\tilde{u}_t(x): (t,x) \in A)$ are mutually absolutely continuous.
\end{cor}
\textbf{Proof} We may suppose that the initial functions $u_0=f$ and 
$\tilde{u}_0=g$ are fixed elements of $\mathcal{E}_{\exp}$, and that the 
two solutions are defined on the same probability space and with respect to 
the same noise $W$.  The case where $u_0$ and $ \tilde{u}_0$ are random 
then follows by using the Markov property at time zero.  We may also suppose 
that $A$ is a rectangle and choose a $\mathbf{C}^{\infty}$ function 
$\psi_t(x)$ that equals $1$ on $A$ and has compact support inside
$(0,\infty) \times \mathbf{R}$.  Define
\[
v_t(x) = u_t(x) + \psi_t(x) \int G_t(x-y) (g(y)-f(y)) dy.
\]
Then using the representation (\ref{solution}) we see that 
$v_t(x) = \tilde{u}_t(x)$ for $(t,x) \in A$.  Also $v_0 =f$ and it is 
straightforward to check that $(v_t(x))$ is a solution to 
(\ref{perturbation}) with
\[
h_t(x) = \left(\frac{\partial}{\partial t} 
- \frac{\partial^2}{\partial^2x} \right) \left(\psi_t(x) 
        \int G_t(x-y) (g(y)-f(y)) dy\right).
\]
Note that $h_t(x)$ is smooth, deterministic  and of compact support and so 
certainly satisfies the hypothesis of Lemma (\ref{girsanov}).  The result 
then follows from Lemma (\ref{girsanov}) by taking $T$ large enough that 
$A \subseteq [0,T] \times \mathbf{R}$.  \qed
\begin{cor} \label{absolutecorollary2}
Suppose $(u_t(x))$ is a solution to (\ref{spde}) and $z \in \mathbf{R}^d$.
For any compact set $A \subseteq (0,\infty) \times \mathbf{R}$ the laws of the 
fields $(u_t(x): (t,x) \in A)$ and $(z + u_t(x): (t,x) \in A)$ are mutually 
absolutely continuous.
\end{cor}
\textbf{Proof} The proof is similar to the proof of the previous corollary,
but one defines $v_t(x) = u_t(x) + z \, \psi_t(x)$ and changes $h_t(x)$ 
accordingly.  \qed

\vspace{.1in}

Our next aim is to show sufficient absolute continuity to 
allow us to transfer our results from the random string to the random loop. 
A continuous adapted  process 
$(\tilde{u}_t(x): t \geq 0, \,  x \in \mathbf{T})$ is a solution to the 
random loop form of (\ref{spde}) if it satisfies (\ref{solution}) where 
$G_t(x)$ is replaced by the Green's function for the heat equation on the 
circle and the stochastic integral is only over the circle $\mathbf{T}$. 
This requires only a white noise $W(dx\,dt)$ defined on 
$t \geq 0, x \in \mathbf{T}$. 

The corollary below implies that the properties we prove about the random 
string in theorems 1,2 and 3 hold also for the random loop. 
\begin{cor} \label{absolutecorollary4}
Suppose $(u_t(x):t \geq 0, \; x \in \mathbf{R})$ is a solutions to 
(\ref{spde}) and 
$(\tilde{u}_t(x):t \geq 0, x \in \mathbf{T})$ is a solution to 
(\ref{spde}) on the circle.  For any compact set
$A \subseteq (0,\infty) \times (0,1)$ the laws of the fields 
$(u_t(x): (t,x) \in A)$ and
$(\tilde{u}_t(x): (t,x) \in A)$ are mutually absolutely continuous.
\end{cor}
\textbf{Proof} We may suppose that the initial functions 
$u_0=f \in \mathcal{E}_{exp}$ and $\tilde{u}_0=g \in C(\mathbf{T})$ are 
deterministic.  The case where $u_0$ and $\tilde{u}_0$ are random then 
follows by using the Markov property at time zero.  We also suppose that 
they are defined on the same probability space and the noise driving 
$(\tilde{u}_t(x))$ is the restriction to the circle of the 
noise $W$ driving $(u_t(x))$.

We use a standard symmetry trick to extend the solution $(\tilde{u}_t(x))$ 
over the real line.  We may extend the solution to  
$(\tilde{u}^{(per)}_t(x):t \geq 0, \, x \in \mathbf{R}) $ by making it 
periodic with period one.  We also extend the noise to a noise 
$W^{(per)}(dx\,dt)$ over the whole line by making it periodic.  Note that  
$\tilde{u}^{(per)}_t(x) = \tilde{u}_t(x)$ and $W^{(per)}(dx\,dt)=W(dx\,dt)$ 
for 
$t \geq 0, \; x \in \mathbf{T}$.  Then $(\tilde{u}^{(per)}_t(x)) $ 
satisfies (\ref{solution}) over the whole line, with the Green's function 
for the whole line but with the periodic noise $W^{(per)}(dx\,dt)$.

We again take a $\mathcal{C}^{\infty}$ function $\psi_t(x)$ that equals $1$ 
on $A$ and still has compact support inside $(0,\infty) \times (0,1)$.  Define
\begin{eqnarray*}
v_t (x) & = & u_t(x) + \psi_t(x) \int G_t(x-y) 
\left(g^{(per)}(y)-f(y) dy \right) \\
& & \hspace{.1in}  + \psi_t(x) \int_{0}^{t} \int
        G_{t-s}(x-y) \left( W^{(per)}(dy\,ds)-W(dy\,ds) \right).
\end{eqnarray*}
Then using the representation (\ref{solution}) we see that 
$v_t(x) = \tilde{u}_t(x)$ for $(t,x) \in A$.  Also 
$v_0 =f$ and it is straightforward to check that $(v_t(x))$ is a solution 
to (\ref{perturbation}) with
\begin{eqnarray*}
h_t(x) &=& \left(\frac{\partial}{\partial t} 
- \frac{\partial^2}{\partial^2x} \right)  
\left( \psi_t(x) \int G_t(x-y) \left( g^{(per)}(y)-f(y) \right) dy  \right) \\
& & \hspace{.2in}   + \left( \frac{\partial}{\partial t} 
- \frac{\partial^2}{\partial^2x} \right)  
\bigg( \psi_t(x) \int_{0}^{t} \int
 G_{t-s}(x-y)   \\
&& \hspace{1.5in} \cdot\left( W^{(per)}(dy\,ds)-W(dy\,ds) \right) \bigg).
\end{eqnarray*}
Note that $h_t(x)$ has compact support.  We claim that $h_t(x)$ is also 
smooth.  The only term in $h_t(x)$ for which this is not clear is 
the stochastic integral 
\[ 
I(t,x) = \int_{0}^{t} \int  G_{t-s}(x-y) \left( W^{(per)}(dy\,ds)-W(dy\,ds) 
\right).
\]
However since $ W^{(per)}(dy\,ds)-W(dy\,ds)=0$ for $y \in (0,1)$ the 
function $I(t,x)$ solves the 
deterministic heat equation in the region $[0,\infty) \times (0,1)$, with 
zero initial conditions
and continuous random boundary values $I(t,1)$ and $I(t,0)$.  Hence it 
is smooth in this region and since $\psi_t(x)$ is also supported in this 
region the claim follows.

Since $ W^{(per)}(dy\,ds)-W(dy\,ds)=0$ for $y \in (0,1)$, the perturbation 
$(h_t(x))$ is adapted to the $\sigma$-field $\mathcal{G}$ generated by the 
noise $W(f)$ for $f$ supported outside $(0,1) \times [0,\infty)$.  Hence the  
integrand $h_t(x)$ is independent of the noise 
$(W(dx\,dt): t \geq 0,\, x \in \mathbf{T})$,
and we can apply part (b) of Lemma \ref{girsanov}.  \qed
\begin{cor} \label{absolutecorollary3}
Suppose $(u_t(x))$ is a solution to (\ref{spde}).  Suppose also that  $A^+$ 
is a compact set in the half space
$\mathbf{H}^+ = (0,\infty) \times (0,\infty)$ and $A^-$ is a compact set in 
the half space $\mathbf{H}^- = (0,\infty) \times (-\infty,0)$.  Then the 
law of the pair of fields 
\[
\left((u_t(x): (t,x) \in A^+), (u_t(x): (t,x) \in A^-)\right)
\]
is mutually absolutely continuous with respect to the law of  
\[
\left((U_t(x): (t,x) \in A^+), (\tilde{U}_t(x): (t,x) \in A^-)\right)
\]
where $(U_t(x))$ and $(\tilde{U}_t(x))$ are independent copies of the 
stationary pinned string.
\end{cor}
\textbf{Proof} We may suppose that the initial function $u_0 =f$ is 
deterministic.  Suppose also $(u_t(x))$ is driven by a noise $W(dx\,dt)$.  
On the same probability space, construct solutions $(u^{+}_t(x))$ 
(respectively $(u^{-}_t(x))$) to (\ref{spde}) on the half space $\mathbf{H}^+$
(respectively $\mathbf{H}^-$) with zero initial conditions and with Dirichlet 
boundary conditions along the axis $\{x=0\}$.  The noise driving $(u^{+}_t(x))$ 
(respectively $(u^{-}_t(x))$) is $(W(dx\,dt): (t,x) \in \mathbf{H}^+)$ (respectively 
$(W(dx\,dt): (t,x) \in \mathbf{H}^-)$).  We can represent the solution $u^{+}_t(x)$ by
\[
u_t^{+}(x) = \int^t_0 \int G_{t-s}(x-y) W^{+}(dy\,ds), \quad \mbox{for 
$t\geq0, \; x\geq 0$}
\]
where $W^{+}(dx\,dt)$ is the odd extension of the noise 
$(W(dx\,dt): (t,x) \in \mathbf{H}^+)$
defined by 
\[
W^{+}([-b,-a] \times [s,t]) = - W([a,b] \times [s,t]), 
\quad \mbox{for all $0 \leq s \leq t,\, 0 \leq  a \leq b$.}
\]
A similar representation holds for $u^{-}_t(x)$ using an odd extension of  
$(W(dx\,dt): (t,x) \in \mathbf{H}^- )$.  Note that $(u^{+}_t(x))$ and $(u^{-}_t(x))$ 
are independent. 

Now choose  $\psi_t^{+}(x)$ (respectively $\psi^{-}_t(x)$)  smooth, equal to
$1$ on $A^+$ (respectively
on $A^-$), and  supported in $\mathbf{H}^+$ (respectively in $\mathbf{H}^-$).  Define 
\begin{eqnarray*}
v_t(x) & = & u_t(x) + \psi^+_t(x) \bigg( - \int G_t(x-y) f(y) dy  \\
&&\qquad 
+ \int^t \! \int G_{t-s}(x-y) \left( W^+(dy\,ds) - W(dy\,ds) \right) \bigg) \\
&& + \psi^-_t(x) \bigg( - \int G_t(x-y) f(y) dy   \\
&&\qquad + \int^t \! \int G_{t-s}(x-y) 
\left( W^-(dy\,ds) - W(dy\,ds) \right) \bigg).
\end{eqnarray*}
We now argue as in Corollary \ref{absolutecorollary4}.  $v_t(x)$ agrees with 
$u^+_t(x)$ on $A^+$ and with
$u^-_t(x)$ on $A^-$.  Also it solves (\ref{perturbation}) with a suitable 
drift $h_t(x)$ that satisfies
the assumptions of Lemma \ref{girsanov} (b).  So the law of the pair 
\[
\left((u_t(x): (t,x) \in A^+), (u_t(x): (t,x) \in A^-)\right)
\]
is mutually absolutely continuous with respect the law of the pair 
\[
\left((u^+_t(x): (t,x) \in A^+), (u^-_t(x): (t,x) \in A^-)\right).  
\]
But this second pair is independent and a similar argument to the above 
shows that $(u^+_t(x): (t,x) \in A^+)$ (respectively 
$(u^-_t(x): (t,x) \in A^-)$) is absolutely continuous with respect to 
$(U_t(x): (t,x) \in A^+)$ (respectively $ (\tilde{U}_t(x): (t,x) \in A^+)$).  
\qed
%
%%%%%%%%%%%%%%%%%%%%%%%%%%%%%%%%%%%%%%%%%%%%%%%%%%%%
\section{Hitting Points} \label{s4}
\setcounter{equation}{0}
%
%%%%%%%%%%%%%%%%%%%%%%%%%%%%%%%%%%%%%%%%%%%%%%%%%%%%
For a $R^d$-valued function $u_t(x)$ indexed over 
$(t,x) \in A \subseteq [0,\infty) \times \mathbf{R} $, we say  that 
$(u_t(x): (t,x) \in A)$ hits the point $z \in \mathbf{R}^d$ if $u_t(x) =z$ 
for some $(t,x) \in  A$. 
The aim of this section is to prove the following result. 
\begin{theorem} \label{hitting1}
Suppose $(u_t(x): t \geq 0, x \in \mathbf{R})$ is a solution to (\ref{spde}). 
\begin{itemize}
\item[(a)] If $d \leq 5$ then 
$ P \left( (u_t(x): (t,x) \in A) \; \mbox{hits $z$}\right) >0 $ for all 
$z \in \mathbf{R}^d$ and all $A \subseteq [0,\infty) \times \mathbf{R} $ 
with non-empty interior.
\item[(b)] If $d \geq 6$ then 
$ P \left( (u_t(x): t>0, x \in \mathbf{R}) \; \mbox{hits $z$} \right) = 0$
for all $z \in \mathbf{R}^d$.
\end{itemize}
\end{theorem}
To prove the theorem we shall need the following lemma which gives 
covariance estimates on the events of the stationary pinned string hitting 
a small ball. 
\begin{lemma} \label{covarianceestimates1}
There exist constants $0<c_2, c_3<\infty$, depending only on the dimension 
$d$,  so that the following 
bounds hold: for all $s,t \in [1,2]$, $x,y \in [-2,2]$ and $\delta \in (0,1]$
\begin{eqnarray} \label{lemma-assertion-1} 
&& P\Big(U_t(x) \in B_{\delta}(0) \Big)  \geq  c_2 \delta^d,\\
\label{lemma-assertion-2}
&& P\Big(U_t(x) \in B_{\delta}(0) ,U_s(y) \in B_{\delta}(0) \Big)  \leq  
 c_3 \delta^{2d} \left(|t-s|^{1/2}+|x-y| \right)^{-d/2}.   \nonumber\\
&&
\end{eqnarray}
\end{lemma}
\textbf{Proof of Lemma \ref{covarianceestimates1}.} The Gaussian variable 
$U^{(i)}_t(x)$ has mean zero and variance 
$t^{1/2} F(|x|t^{-1/2}) \geq (2 \pi)^{-1/2}$ for 
$x \in [-2,2],\, t \in [1,2]$.  So it has a density which is bounded above 
by $(2 \pi)^{-1/4} $ and (\ref{lemma-assertion-1}) follows by the 
independence of the coordinates $U^{(i)}_t$.  An analogous upper bound also 
holds.

To prove (\ref{lemma-assertion-2}) we consider the mean zero Gaussian vector
 $(X,Y) = (U^{(i)}_t(x), U^{(i)}_s(y) ) $.  The covariance 
 (\ref{generalcovariance}) gives an expression
for $\sigma_X^2= E[X^2]$, $ \sigma_Y^2 = E[Y^2]$ and
$\rho_{X,Y}^2 = E[(X-Y)^2]$.  The law of $X-E[X|Y]$ is Gaussian and a 
routine calculation shows
its variance is given by   
\begin{equation} 
\mbox{Var}\left(X -E[X|Y]\right) 
=  \frac{\left(\rho_{X,Y}^2 - (\sigma_X - \sigma_Y)^2 \right)
\left((\sigma_X + \sigma_Y)^2 - \rho_{X,Y}^2 \right)}{4 \sigma_X^2}.  
\label{conditionalvariance}
\end{equation}
For mean zero Gaussian $Z$ the probability $P(\mu + Z \in B_{\delta}(0))$ 
is maximized at $\mu = 0$. 
So we can bound the probability  
\[
P \Big( X \in B_{\delta}(0) | Y \Big)\le
 C \delta \cdot\Big(\mbox{Var}(X -E[X|Y] ) \Big)^{-1/2}
\]
and hence obtain
\begin{equation} \label{estimate103}
P\Big(X \in B_{\delta}(0) ,Y \in B_{\delta}(0) \Big) 
\leq C \delta^2 \Big(\mbox{Var}(X -E[X|Y] ) \Big)^{-1/2}.
\end{equation}
The covariance (\ref{generalcovariance}) implies that $\sigma_X$ is bounded 
and bounded away from zero, for $t \in [1,2]$ and $|x| \leq 2$.  The 
inequality (\ref{boundoncovariance}) implies that
$ \rho_{X,Y}^2 \geq c_1 \left( |t-s|^{1/2} + |x-y| \right).$ The 
differentiability of $F(z)$ and the mean value theorem combine to show that  
$ | \sigma_X - \sigma_Y| \leq C \left( |t-s| + |x-y| \right).  $
Using these bounds in (\ref{conditionalvariance}) shows there exists $C>0$ 
and $\epsilon>0$ so that
\begin{equation}
\label{variance113}
\mbox{Var}\Big(X -E[X|Y]\Big) \geq C \left( |t-s|^{1/2} + |x-y| \right) 
\end{equation}
whenever $t,s \in [1,2]$, $x,y \in [-2,2]$ and $|t-s| + |x-y| \leq \epsilon$.
The variance $\mbox{Var}(X-E(X|Y))$ is a continuous function of 
$s,t \in [1,2]$ and $x,y \in [-2,2]$. 
It vanishes in this region only on $s=t,x=y$ and hence is bounded below 
when $|t-s| + |x-y| \geq \epsilon$.  So, changing the constant $C$ if 
necessary, the lower bound (\ref{variance113}) holds without the 
restriction  $|t-s| + |x-y| \leq \epsilon$.  Substituting 
(\ref{variance113}) into (\ref{estimate103}), and using the independence of 
coordinates $U^{(i)}$ gives the desired bound.  \qed

\vspace{.1in}

\noindent \textbf{Proof of Theorem \ref{hitting1} in $d \leq 5$.}
We start with a series of five easy reductions.  First, the projection of a 
solution into a lower dimension is still a solution.  Thus we need only 
argue in dimension $d=5$.  Second, it is enough to prove the result when $A$ 
is a compact rectangle in $(0,\infty) \times \mathbf{R}$.  Third, by the 
absolute continuity from Corollary \ref{absolutecorollary1} it is enough to 
prove the result for the stationary pinned string.  Fourth, the absolute 
continuity from Corollary \ref{absolutecorollary2} shows that 
$ P \left( (u_t(x): (t,x) \in A) \; \mbox{hits $z$}\right) $ is either zero 
for all $z$ or strictly positive for all $z$.  So it is enough to prove the 
result when $z=0$.  Fifth and finally, the scaling of the stationary pinned 
string implies that it is enough to consider the rectangle 
$A= [1,2] \times [0,1]$.  To see this, note that if  
$ P \left( (U_t(x): (t,x) \in A) \; \mbox{hits $0$}\right) >0 $ then, as 
above, the absolute continuity results imply for a solution $(u_t(x))$ to 
(\ref{spde}) started at any  $f \in \mathcal{E}_{exp}$, and for 
any $z \in \mathbf{R}^d$, that  
$ P \left( (u_t(x): (t,x) \in A) \; \mbox{hits $z$}\right) >0 $.  Then by 
applying the Markov property at a time $t_0$ one sees that 
\[
P \Big( (U_{t+t_0}(x+x_0): (t,x) \in A) \; \mbox{hits $0$}\Big) >0 
\]
for any $t_0 \geq 1$ and $x_0 \in \mathbf{R}$.  The scaling of the 
stationary pinned string gives, for any $0<r<s$ and $a<b$,
\begin{eqnarray*}
\lefteqn{P \Big( (U_t(x): (t,x) \in [r,s] \times [a,b]) 
	\; \mbox{hits $0$}\Big)    }\\
&=& P \left( (U_t(x): (t,x) \in [L^4 r, L^4 s] \times [L^2 a, L^2 b]) \; 
\mbox{hits $0$}\right) \\
& \geq & P \left( (U_{t+(L^4 r-1)}(x + L^2 a): (t,x) \in A) \; 
\mbox{hits $0$}\right) \\
&>&0
\end{eqnarray*}
provided we pick, as we may, $L$ large enough that $L^4 r \geq 1$, 
$L^4 (s-r) \geq 1$, and $L^2 (b-1) \geq 1$.

Now fix $A= [1,2] \times [0,1]$ for the rest of this proof.  Recall the grid 
of points $t_{i,n}$ and $x_{i,n}$ defined in (\ref{grid}).  Define events 
\[
\mathcal{B}_{i,j,n} = \left\{ U_{1+t_{i,n}}(x_{j,n}) \in B_{\delta_n}(0) 
\right\}, \quad  
\mathcal{B}_n=\bigcup_{i=1}^{2^{4n}} \bigcup_{j=1}^{2^{2n}} 
\mathcal{B}_{i,j,n}.  
\]
We shall show that $P( \mathcal{B}_n ) \geq p_0>0$ for all $n$.  Then, using 
continuity of $U$ and the compactness of $A$,
\[ 
P \Big((U_t(x): (t,x) \in A) \; \mbox{hits the point $0$} \Big) \geq 
P(\mathcal{B}_n \; \mbox{infinitely often}) \geq p_0.
\]
We shall apply Lemma \ref{lower-bound} to the events $\mathcal{B}_{i,j,n}$. 
First, 
(\ref{lemma-assertion-1}) applied in dimension $d=5$ implies that 
\begin{equation} \label{lower-1}
\sum_{i=1}^{2^{4n}}\sum_{j=1}^{2^{2n}} P\left(\mathcal{B}_{i,j,n} \right) 
\geq  
c_2 2^{6n} \delta_n^5 = c_2. 
\end{equation}
Second, using  (\ref{lemma-assertion-2}), 
\begin{eqnarray} \label{inclusion-exclusion-1}
\lefteqn{\sum_{i=1}^{2^{4n}} \sum_{j=1}^{2^{2n}} 
\sum_{\tilde{i}=1}^{2^{4n}} \sum_{\tilde{j}=1}^{2^{2n}} P\left( 
\mathcal{B}_{i,j,n} \cap \mathcal{B}_{\tilde{i}, 
\tilde{j},n}\right) \mathbf{1} \left((i,j) \neq (\tilde{i}, 
\tilde{j}) \right)   }\\
& \leq & 2 \sum_{i=1}^{2^{4n}} \sum_{j=1}^{2^{2n}} \sum_{k=0}^{2^{4n}}
\sum_{\ell=-2^{2n}}^{2^{2n}} 
P\left(\mathcal{B}_{i,j,n}\cap\mathcal{B}_{i+k,j+\ell,n}\right)  
\mathbf{1}\left((k,\ell) \neq (0,0) \right)  \nonumber \\ 
& \leq & 2 c_3 2^{6n}  \delta_n^{10}  \sum_{k=0}^{2^{4n}} 
\sum_{\ell=-2^{2n}}^{2^{2n}}
\left(\left| k2^{-4n}\right|^{1/2} +\left| \ell 2^{-2n}\right| 
\right)^{-5/2} \mathbf{1} \left((k,\ell) \neq (0,0) \right)    \nonumber\\
& \leq & 2^2 c_3 2^{11n}  \delta_n^{10}  \sum_{k=0}^{2^{4n}} 
\sum_{\ell=0}^{2^{2n}}
  \left( k^{1/2}+ | \ell | \right)^{-5/2}  \mathbf{1}\left((k,\ell) \neq (0,0) 
  \right)     \nonumber  \\  
& \leq & 2^2 3^{5/2} c_3 2^{11n}  \delta_n^{10}  \sum_{k=1}^{2^{4n}+1} 
\sum_{\ell=1}^{2^{2n}+1}
 \left( k^{1/2}+| \ell | \right)^{-5/2}  \nonumber  \\
& \leq & 2^2 3^{5/2} c_3 2^{11n}  \delta_n^{10}  
\int_{0}^{2^{4n}+1}\int_{0}^{2^{2n}+1}
 \left(x^{1/2}+y \right)^{-5/2} dydx     \nonumber \\
& \leq   &  2^3 3^{3/2} c_3 2^{11n} \delta_n^{10}  \int_{0}^{2^{4n}+1} 
x^{-3/4} dx  \nonumber \\
& = & 2^{21/4} 3^{3/2} c_3 2^{12n} \delta_n^{10} \nonumber \\
&  \leq & 2^6 3^{2} c_3 \nonumber
\end{eqnarray}
Using Lemma \ref{lower-bound}, together with (\ref{lower-1}) and 
(\ref{inclusion-exclusion-1}), we obtain
\[
P(\mathcal{B}_{n}) \geq  \frac{c_2^2}{1 + 2^6 3^2 c_3} >0 
\]
for all $n \geq 1$. 
This completes the proof that points can be hit in dimensions $d \leq 5$. 
The reader can check that the  above proof would fail if we replace $d=5$ 
by $d=6$.

\vspace{.1in}

\noindent \textbf{Proof of Theorem \ref{hitting1} in $d \geq 6$.} We again 
make some reductions.  By considering projections of the string into lower 
dimensions, it is enough to consider dimension $d=6$.  It is enough to show 
that $ P(u_t(x) = z \; \mbox{for some $(t,x) \in A$})=0 $ for a bounded 
rectangle $A$.  It is then  enough to consider the stationary pinned string 
and again, using scaling, it is enough to consider
$A=[0,1) \times [0,1)$.  Finally, since the probability 
$ P(U_t(x) = z \; \mbox{for some $(t,x) \in A$})$ is either zero for all 
$z$ or strictly positive for all $z$, the problem can be tackled  by 
studying the range of the process, defined by 
\[
U(A)=\{U_t(x): (t,x) \in A\} \subseteq \mathbf{R}^6.
\]
Indeed, if we denote the Lebesgue measure of $U(A)$ by $m(U(A))$ then
\[
E[m(U(A))] = \int_{\mathbf{R}^6} P(U_t(x) = z \; \mbox{for some 
$(t,x) \in A$}) \; dz,
\]
which is zero if and only if the integrand is identically zero.

Subdivide $A$ into eight disjoint rectangles $A_1, \ldots, A_8$, each a 
translate of $[0,1/4) \times [0,1/2)$.  The scaling and translation 
invariance of the stationary pinned string in dimension $d=6$
implies that $E[m(U(A_i))]= (1/8) E[m(U(A))]$ for all $i=1,\ldots,8$. 
However by an 'inclusion-exclusion' type argument
\[ 
m(U(A)) \leq  \sum_{i=1}^8 
m(U(A_i)) - m \Big( U(A_1) \cap U(A_2) \Big).
\]
Taking expectation of both sides shows that $ E[m(U(A_1) \cap U(A_2))] =0.  $ 
We may suppose that
$ A_1 = [0,1/4) \times [0,1/2)$ and $A_2 = [1/4,1/2) \times [0,1/2)$.  Let 
$\mathcal{H}$ be the $\sigma$-field generated by 
$(U_{1/4}(x):  x \in \mathbf{R})$.  Next, we use the Markov property of 
solutions and time reversal for the stationary pinned string. 
Conditioned on $\mathcal{H}$, the laws of 
$(U_t(x): 1/4 \leq t \leq 1/2, \, x \in \mathbf{R})$ and 
$(U_{1/4-t}(x): 1/4 \leq t \leq 1/2, \, x \in \mathbf{R})$ are identical 
and independent.  So
\begin{eqnarray}
0 &=& E\Big[m(U(A_1)\cap U(A_2))\Big]  \nonumber \\
&=& \int_{\mathbf{R^6}} E\Big[ \mathbf{1}( x \in U(A_1)) \mathbf{1} (x \in 
U(A_2) \Big]dx 
\nonumber  \\   
&=& E \left( \int_{\mathbf{R^6}} E\left[ \mathbf{1}( x \in U(A_1)) 
\mathbf{1}( x \in U(A_2))
 \Big| \mathcal{H} \right]dx \right)  \nonumber \\
&=& E\left(\int_{\mathbf{R^6}} E\left[\mathbf{1}( x \in U(A_1)) 
\Big|\mathcal{H}\right]
 E \left[ \mathbf{1}( x \in U(A_2)) \Big|\mathcal{H}\right]dx\right)   
 \nonumber \\
&=&  E\left(\int_{\mathbf{R^6}}
  E \left[ \mathbf{1}( x \in U(A_1)) \Big|\mathcal{H}\right]^2 dx\right)  
  \label{conditioned}
\end{eqnarray}
This implies that 
$ E\left[\mathbf{1}( x \in U(A_1)) \Big|\mathcal{H}\right]=0 $ for almost 
every $x$, almost surely.  But then we have
\[
E \left[ m(U(A)) \right]=8 E \left[m(U(A_1))\right] = E \left( 
\int_{\mathbf{R}^6} 
E\left[\mathbf{1}( x \in U(A_1)) \Big|\mathcal{H}\right] dx \right) =0
\]
and therefore $m(U(A))=0$ almost surely, which concludes the proof .  \qed
%
%%%%%%%%%%%%%%%%%%%%%%%%%%%%%%%%%%%%%%%%%%%%%%%%%%%%
\section{Double points} \label{s5}
\setcounter{equation}{0}
%
%%%%%%%%%%%%%%%%%%%%%%%%%%%%%%%%%%%%%%%%%%%%%%%%%%%%
We consider two kinds of double points.  For a $\mathbf{R}^d$ valued 
function $u_t(x)$, we say  that 
$(u_t(x): (t,x) \in A)$ has a double point at $z \in \mathbf{R}^d$ if
there exist $(t,x), \, (t,y) \in A$, with $x \neq y$, so that 
$u_t(x)=u_t(y)=z$.
We say that the range of the function $(u_t(x): (t,x) \in A)$ has a double 
point $z$
 if there exist $(t,x), \, (s,y) \in A$, with $(t,x) \neq (s,y)$, such  
 that $u_t(x)=u_s(y)=z$. 
The aim of this section is to prove the following result. 
\begin{theorem} \label{multiple-points} 
Suppose $(u_t(x):t \geq 0, x \in \mathbf{R})$ is a solution to (\ref{spde}),
and let $A \subseteq (0,\infty) \times \mathbf{R}$ have non-empty interior.  
The following statements hold almost surely.  
\begin{itemize}
\item[(a)]  If $d\le 7$, then $(u_t(x): (t,x) \in A)$ has a double point. 
\item[(b)] If $d \geq 8$, then $(u_t(x):t>0, x \in \mathbf{R})$ has no double 
points.  
\item[(c)] If $d \leq 11$, then the range of $ (u_t(x): (t,x) \in A)$ 
has a double point.  
\item[(d)] If $d \geq 12$, then the range of $ (u_t(x):t>0, x \in \mathbf{R})$ 
has no double points.  
\end{itemize}
\end{theorem}
\textbf{Remarks} 

\noindent {\bf 1.} One could also consider double points at a fixed time, 
that is, fix $t>0$ and
ask if there exist $x \neq y$ so that $u_t(x)=u_t(y)$.  However, the 
covariance structure (\ref{spatialcovariance}) implies that the process 
$x \to U_t(x)-U_t(0)$ is a two sided Brownian motion.  
It is well known  that there are double points, with non-zero 
probability, if and only if $d <4 $.  Absolute continuity then shows the 
same holds true for general solutions to (\ref{spde}).

\noindent {\bf 2.}  Parts (a) and (c) on the existence of double points 
follow by an inclusion-exclusion argument similar to that in theorem 
\ref{hitting1}.  We illustrate this by giving the argument for part (c), 
which is the more complicated, and leave the details of part (a) to the 
reader.  In proof of non-existence, we need some small tricks to reduce 
the argument to the scaling property of the stationary string.

\vspace{.1in} 

\noindent \textbf{Proof of Theorem \ref{multiple-points} (c): existence of 
double points of the range  in  dimensions $d \leq 11$.} We can again make 
various reductions, arguing as in the proof of theorem \ref{hitting1}.  By 
projection it is enough to argue in dimensions $d = 11$.  It is enough to 
consider bounded $A$, and hence by absolute continuity, enough to consider 
the stationary pinned string.  Scaling and translation invariance for the 
stationary string again imply it is enough to consider one fixed rectangle, 
say $A=[0,4] \times [0,1]$. 

For the rest of this proof we set
\[
A_1= [0,1] \times [0,1] \quad \mbox{and} \quad A_2 = [3,4] \times [0,1]
\] 
and $\delta_n=2^{-12n/11}$.  Define the events 
\[
\mathcal{B}_{i,j,k,\ell,n}= \left\{ 
U_{t_{i,n}}(x_{j,n})-U_{3+t_{k,n}}(x_{\ell,n}) \in B_{\delta_n}(0) 
\right\}, \quad  
\mathcal{B}_n=\bigcup_{i,k=1}^{2^{4n}}\bigcup_{j,\ell=1}^{2^{2n}}
\mathcal{B}_{i,j,k,\ell,n}. 
\]
We will show $P(\mathcal{B}_n) \geq p_0 >0$ for all $n$.  Then,
by continuity and compactness, we have 
\begin{eqnarray*}
\lefteqn{P \left( \mbox{the range of 
$ (U_t(x): (t,x) \in [0,4] \times [0,1]) $ has a double point}\right) } \\
& \geq & P( \left\{U_t(x): (t,x) \in A_1 \right\} \cap \left\{U_{t}(x): 
(t,x) \in A_2 \right\} \neq \emptyset)\\ & \geq & P (\mathcal{B}_n \, 
\mbox{infinitely often}) \geq p_0.
\end{eqnarray*}
We need the following lemma on the covariance structure of the events 
$\mathcal{B}_{i,j,k,\ell,n}$. 
\begin{lemma} \label{covarianceestimates2}
Suppose that $s_i,t_i,x_i,y_i\in [0,1]$ for $i=1,2$.  There exist constants 
$0 < c_4, c_5 < \infty$, depending only on the dimension $d$, so
 that for all $0<\delta\le 1$
\begin{eqnarray}
\label{ball-1}
\lefteqn{P\left(U_{t_1}(x_1)-U_{3+s_1}(y_1) \in B_{\delta}(0) \right) \geq 
c_4 \delta^d}    \\
\lefteqn{P\left( U_{t_1}(x_1)- U_{3+s_1}(y_1)\in B_{\delta}(0), 
U_{t_2}(x_2) - U_{3+s_2}(y_2) 
\in B_{\delta}(0)  \right)}  \nonumber \\
 & \leq & c_5 \delta^{2d} \left( |t_1-t_2|^{1/2}+|s_1-s_2|^{1/2}+|x_1-x_2| 
 + |y_1-y_2| \right)^{-d/2}.
  \label{ball-2}
\end{eqnarray}
\end{lemma}
We delay the proof of this lemma until after we complete the main argument. 
Using estimate (\ref{ball-1}), we conclude that
\begin{equation} \label{first-bound}
\sum_{i,k=1}^{2^{4n}}\sum_{j,\ell=1}^{2^{2n}} 
P\left(\mathcal{B}_{i,j,k,\ell,n}\right)
         \geq  c_4 2^{12n} \delta_n^{11}  = c_4.
\end{equation}
Using estimate (\ref{ball-2}), we find that
\begin{eqnarray*}
\lefteqn{ \sum_{i_1, i_2, k_1, k_2=1}^{2^{4n}}  \sum_{j_1, j_2, \ell_1, 
\ell_2=1}^{2^{2n}} 
P \left(\mathcal{B}_{i_1,j_1,k_1,\ell_1,n} \cap 
\mathcal{B}_{i_2,j_2,k_2,\ell_2,n} \right) 
\mathbf{1}_{( (i_1,j_1,k_1,\ell_1) \neq  
(i_2,j_2,k_2,\ell_2) )}  }  \\
& \leq & 2 \sum_{i_1,k_1=1}^{2^{4n}}  \sum_{j_1, \ell_1=1}^{2^{2n}} 
\sum_{i_2,k_2=-2^{4n}}^{2^{4n}}   \sum_{j_2,\ell_2=-2^{2n}}^{2^{2n}}
 P \left( \mathcal{B}_{i_1,j_1,k_1,\ell_1,n} \cap  
\mathcal{B}_{i_1+i_2,j_1+j_2,k_1+k_2,\ell_1+ \ell_2,n} \right)  \\
&& \hspace{2.5in} \cdot\mathbf{1}_{( (i_2,j_2,k_2,\ell_2) \neq (0,0,0,0) )} \\
& \leq & c_5 2^{n+1} \delta_n^{22} \sum_{i_2,k_2=-2^{4n}}^{2^{4n}}  
\sum_{j_2,\ell_2=-2^{2n}}^{2^{2n}} 
\left( |i_2|^{1/2} + |j_2| + |k_2|^{1/2} + |\ell_2| \right)^{-11/2}  \\
&& \hspace{2.5in} \cdot\mathbf{1}_{( (i_2,j_2,k_2,\ell_2)  \neq (0,0,0,0) )}.
\end{eqnarray*}
It is straightforward, as in the proof of theorem \ref{hitting1}, to bound 
this quadruple sum by 
a constant, independent of $n$.  Using  this and (\ref{first-bound}) in 
Lemma \ref{lower-bound} completes the proof that 
$P(\mathcal{B}_n) \geq p_0>0$. 

One way to show that the probability of double points in the range is 
actually one is to use scaling and a zero one law.  Alternatively one can 
use the following argument. 
\begin{eqnarray*}
\lefteqn{ P \Big( \mbox{ $ (U_t(x): (t,x) \in [0,4] \times [0,1]) $ has no 
double points}\Big)} \\
&=& P \Big( \mbox{ $ (U_t(x): (t,x) \in [0,16] \times [0,2]) $ has no 
double points}\Big)
\quad \mbox{(by scaling)} \\
& \leq  & P \Big( \mbox{ $ (U_t(x): (t,x) \in [0,4] \times [0,1]) $ } \\
&& \qquad \mbox{ and
$(U_t(x): (t,x) \in [12,16] \times [0,1]) $ have no double points}\Big) \\
& < & P \Big( \mbox{ $ (U_t(x): (t,x) \in [0,4] \times [0,1]) $ has no 
double points}\Big)
\end{eqnarray*} 
The strict inequality in the last line follows by applying the Markov 
property at time $t=4$, the absolute continuity results and the translation 
invariance of the stationary string. 

\vspace{.1in}

\noindent \textbf{Proof of Lemma \ref{covarianceestimates2}.}
The proof follows the argument used for Lemma \ref{covarianceestimates1}, 
with the change that we now consider the Gaussian pair
\[
(X,Y) = \left(U_{t_1}(x_1)- U_{3+s_1}(y_1), U_{t_2}(x_2) - U_{3+s_2}(y_2) 
\right).
\]
The covariance (\ref{generalcovariance}) implies that $\sigma_X$ and 
$\sigma_Y$ are bounded above and away from zero as $s_i,t_i,x_i,y_i$ range 
over $[0,1]$.  Using the identity
\[
(a-b+c-d)^2 = (a-b)^2 + (c-d)^2 + (a-d)^2 +(b-c)^2 -(a-c)^2 -(b-d)^2
\]
we may use (\ref{generalcovariance}) to find, for $t_1 \neq t_2$ and 
$s_1 \neq s_2$,
\begin{eqnarray}
\label{parallelogram}
\rho_{X,Y}^2 & = & E \left( \left(U_{t_1}(x_1)-  U_{t_2}(x_2) +  
U_{3+s_2}(y_2) - U_{3+s_1}(y_1) \right)^2 \right) \nonumber \\
& = & |t_2-t_1|^{1/2} F \left( |x_2 -x_1| \, |t_2 - t_1|^{-1/2} \right) \\
&& \hspace{.2in} + 
|s_2-s_1|^{1/2} F \left( |y_2 -y_1| \, |s_2 - s_1|^{-1/2} \right) \nonumber \\
& &  \hspace{.2in} + H_{t_1-s_1}(x_1-y_1) + H_{t_2-s_2}(x_2-y_2) \nonumber\\
&& \hspace{.2in} - H_{t_1-s_2}(x_1-y_2) - H_{t_2-s_1}(x_2-y_1)   \nonumber
\end{eqnarray}
where $H_r(z) = |3+r|^{1/2} F \left( |z| \cdot |3+r|^{-1/2} \right)$.  Small 
changes are needed
for the cases where $t_1=t_2$ or $s_1=s_2$, but these are easy and left to 
the reader.
The function  $H_r(z)$ is smooth for $r,z \in [-1,1]$.  The last four terms 
on the right hand side 
of (\ref{parallelogram}) are differences of $H$ at the four vertices of a 
parallelogram.  Using the 
mean value theorem twice, these can be expressed as a double integral of  
second derivatives of $H$ over the parallelogram.  Hence the contribution 
of these last four terms is bounded by the size of the second derivatives 
and the area of the parallelogram and is thus at most 
$C( |t_2-t_1|^2 + |s_2-s_1|^2 + |x_2-x_1|^2  
+ |y_2-y_1|^2)$. 
Using (\ref{boundoncovariance}) to bound the first two terms on the right 
hand side of (\ref{parallelogram}) from below we find there exists 
$\epsilon>0$ so that 
\[
\rho_{X,Y}^2 \geq \frac{c_1}{2} \left( |t_2-t_1|^{1/2} + |s_2-s_1|^{1/2} + 
|x_2-x_1| + |y_2-y_1| \right),
\]
whenever $t_i,s_i,x_i,y_i \in [0,1]$ and 
$ |t_2-t_1| + |s_2-s_1| + |x_2-x_1| + |y_2-y_1| \leq \epsilon$.
The rest of the argument exactly parallels that of Lemma 
\ref{covarianceestimates1} and is omitted.  \qed

\vspace{.1in}

\noindent \textbf{Proof of Theorem \ref{multiple-points} (b): non-existence 
of double points  in dimensions $d \geq 8$.} By a projection argument we 
need work only in dimension $d=8$.  It is enough to show there are no double 
points for $(u_t(x): (t,x) \in A) $ for compact 
$A \subseteq (0,\infty) \times \mathbf{R}$, and hence by absolute 
continuity we can work with the stationary pinned string.  We shall show that
\begin{equation} \label{d=8key}
P \left( 0 \in \left\{ U_t(x) -U_t(-y):  t, \, x, \, y \in [1,2) \right\} 
\right) = 0.
\end{equation}
By scaling and translation invariance this implies, for all $t_0,L \geq 0$ 
and $x_0 \in \mathbf{R}$, that
\[
P \left( 0 \in \left\{ U_{t}(x) -U_{t} ( - y):  t \in [t_0,t_0+L^4), \; x,y 
\in [x_0+L^2,x_0+2L^2) 
\right\} \right) = 0.
\]
Taking a countable union of such events shows that there are no double 
points.  Define 
\[ 
V(t,x,y) = U_{1+t}(1+x) - U_{1+t}(-1-y), \quad \mbox{for $t,x,y \in [0,1)$.}
\]
We must show that $P(V(t,x,y)=0 \; \mbox{for some $(t,x,y) \in [0,1)^3$})=0$. 
Define, using an independent copy $\tilde{U}_t(x)$ of the stationary string,
\[
\tilde{V}(t,x,y) = U_{1+t}(1+x) - \tilde{U}_{1+t}(-1-y), \quad \mbox{for 
$t,x,y \in [0,1)$.}
\]
Corollary \ref{absolutecorollary3} implies that the laws of  
$(V(t,x,y): (t,x,y) \in [0,1)^3)$ and 
$(\tilde{V}(t,x,y): (t,x,y) \in [0,1)^3)$ are mutually absolutely 
continuous.  Hence we may work with $\tilde{V}$ in place of $V$.  The 
absolute continuity from Corollary \ref{absolutecorollary2} implies that 
$P(\tilde{V}(t,x,y)=z \; \mbox{for some $(t,x,y) \in [0,1)^3$})$ is either 
zero for all $z$ or strictly positive for all $z$.  Hence, as in theorem 
\ref{hitting1} part (b), it is enough to show that 
$E(m(\tilde{V}([0,1)^3)))=0$.  We can now apply scaling.  Subdivide 
$A=[0,1)^3$ into a disjoint union of 16 rectangles $(A_i: 1=1,\ldots,16)$ 
each of the form 
$A_i = (t_i,x_i,y_i) + A_0 $, with $(t_i,x_i,y_i) \in A$ and 
$A_0=[0,1/4) \times [0,1/2)^2$.  Using the independence of $U$ and 
$\tilde{U}$ and the scaling for the stationary pinned string with 
$L = 2^{-1/2}$, we obtain the following equality in law: 
\begin{eqnarray*}
\lefteqn{m(\tilde{V}(A_i))} \\
& = & m\left( \left\{U_{1+t}(1+x) - \tilde{U}_{1+t}(-1-y) : (t,x,y) \in 
(t_i,x_i,y_i) + A_0 \right\} \right) \\
& \stackrel{\mathcal{L}}{=}  & m\Big( \big\{2^{-1/2}\left(U_{t}(x) - 
\tilde{U}_{s}(-y)\right) :   \\
&& \hspace{1in}   (t,x,y) \in (4+4 t_i, 2+ 2x_i, 2+ 2y_i) + A) 
\big\} \Big) \\
& = & \frac{1}{16} m \left( \left\{U_{t}(x) - \tilde{U}_{s}(-y) : (t,x,y) 
\in (4+4 t_i, 2+ 2x_i, 2+ 2y_i) + A) \right\} \right) \\
& = & \frac{1}{16}  m\Big( \Big\{\Big[ U_{4+4t_i+t}(2x_i+x) - 
U_{3+4t_i}(2x_i)\Big]     \\
&& \qquad\qquad - \Big[ U_{4+4t_i+t}(2y_i+y) -U_{3+4t_i}(2y_i) 
\Big] : (t,x,y) \in  A) \Big\} \Big) \\
& \stackrel{\mathcal{L}}{=}  & \frac{1}{16}  m(\tilde{V}(A)). 
\end{eqnarray*}
The third equality uses the scale factor $(\sqrt{2})^8=16$; the fourth equality 
uses the fact that Lebesgue measure is unchanged by translation; the 
final equality in law uses the translation invariance of the stationary 
pinned string. 

Using the inclusion-exclusion argument from the proof of theorem 
\ref{hitting1}, we obtain 
\[
E \left[ m \left(\tilde{V}(A_i) \cap \tilde{V}(A_j) \right) \right] =0  
\mbox{ for $i \neq j$}.  
\]
The rest of the argument is similar to the theorem \ref{hitting1} part (b). 
We may assume that
\[ 
A_1 = [0,1/4) \times [0,1/2)^2, \quad A_2 = [1/4,1/2) \times [0,1/2)^2.
\] 
Define, for $(t,x,y) \in A_0$, 
\begin{eqnarray*}
V^{(1)}(t,x,y) & = & \left(U_{(5/4) + t}(1+x) -U_{5/4}(0) \right) - 
\left(\tilde{U}_{(5/4) + t}(-1-y)
 - \tilde{U}_{5/4}(0) \right),  \\
V^{(2)}(t,x,y) & = & \left(U_{(5/4) - t}(1+x) -U_{5/4}(0) \right) - 
\left(\tilde{U}_{(5/4) - t}(-1-y)
 -\tilde{U}_{5/4}(0) \right).
\end{eqnarray*}
Note that  
\[ 
m \left(\tilde{V}(A_i) \cap \tilde{V}(A_j) \right)= m \left(V^{(1)}(A_0) 
\cap V^{(2)}(A_0) \right).
\] 
Let $\mathcal{H}$ denote the $\sigma$-field generated by 
$(U_{5/4}(x), \tilde{U}_{5/4}(x): x \in \mathbf{R})$.  Using the Markov 
property, the time reversal and translation invariance of the stationary 
pinned string, the processes $V^{(1)}$ and $V^{(2)}$ are, conditioned on 
$\mathcal{H}$, independent and identically distributed.  
Now we can argue exactly as in (\ref{conditioned}) in the proof of  theorem 
\ref{hitting1} part (b)
to conclude that $E(m(\tilde{V}(A))) = 16 E(m(\tilde{V}(A_1))=0$ which 
finishes the proof
of (\ref{d=8key}).  \qed

\vspace{.1in}

\noindent \textbf{Proof of Theorem \ref{multiple-points} (d): non-existence 
of double points of the range in dimensions $d \geq 12$.}
Only small changes are needed from the proof of part (b).  Again by a 
projection argument we need work only in dimension $d=12$.  It is enough to 
show there are no double points in the range $(u_t(x): (t,x) \in A) $ for 
compact sets $A \subseteq (0,\infty) \times \mathbf{R}$, and hence by absolute 
continuity we can work with the stationary pinned string.  We shall show, 
for any $a \in \mathbf{R}$, that
\begin{equation} \label{d=12key}
P \left( 0 \in \left\{ U_t(x)-U_s(y):  (t,x,s,y) \in [0,1) \times [3,4) 
\times [0,1) \times [a,a+1) 
\right\} \right) = 1.
\end{equation}
By scaling and translation invariance this implies, for all $t_0,L \geq 0$ 
and $x_0,y_0 \in \mathbf{R}$, that
\begin{eqnarray*}
\lefteqn{P \Big( 0 \in \big\{ U_t(x)-U_s(y):  (t,s,x,y) \in 
[t_0,t_0 + L^4)         }\\
&&  \times [t_0 + 3L^4,t_0 + 4L^4) \times [x_0, x_0+L^2) 
\times [y_0,y_0+L^2)  \big\} \Big)   \\
&\hspace{.6in} =&  1.
\end{eqnarray*}
Taking a countable union of such events shows that there are no double 
points  $U_t(x) = U_s(y)$ where $t \neq s$.  Combining this with the result 
of part (b) of the theorem concludes the proof.

Define 
\[ 
V(t,s,x,y) = U_{3+t}(x) - U_{1-s}(a+y), \quad \mbox{for $t,s,x,y \in [0,1)$.}
\]
We must show that 
$P(V(t,s,x,y)=0 \; \mbox{for some $(t,s,x,y) \in [0,1)^4 $})=0$. 
Define, using an independent copy $\tilde{U}_t(x)$ of the stationary string,
\[
\tilde{V}(t,s,x,y) = U_{1+t}(x) - \tilde{U}_{1+s}(y), 
\quad \mbox{for $t,s,x,y \in [0,1)$.}
\]
We claim that the laws of  $(V(t,s,x,y): (t,s,x,y) \in [0,1)^4)$ and 
$(\tilde{V}(t,s,x,y): (t,s,x,y) \in [0,1)^4)$ are mutually absolutely 
continuous.  Indeed, let $\mathcal{H}$ be the $\sigma$-field generated by 
$(U_2(x): x \in \mathbf{R})$.  We can use the Markov property and the time 
reversal property of the stationary pinned string to conclude the following.  
The processes $(U_{3+t}(x): (t,x) \in [0,1)^2)$ and
$(U_{1-s}(a+y): (s,y) \in [0,1)^2)$ are conditionally independent, 
with respect to $\mathcal{H}$.  Also, each is a solution
to (\ref{spde}).  Now the claim follows by applying the absolute continuity 
from Corollary \ref{absolutecorollary1}.

By the claim we may work with $\tilde{V}$ in place of $V$.  The absolute 
continuity from Corollary \ref{absolutecorollary2} implies that 
$P(\tilde{V}(t,s,x,y)=z \; \mbox{for some $(t,s,x,y) \in [0,1)^4$})$ is 
either zero for all $z$ or strictly positive for all $z$.  Hence, as in 
theorem \ref{hitting1} part (b), it is enough to show that 
$E(m(\tilde{V}([0,1)^4)))=0$.  We can now apply scaling.  Subdivide 
$A=[0,1)^4$ into a disjoint union of 64 rectangles $(A_i: 1=1,\ldots,64)$ 
each of the form 
$A_i = (t_i,s_i,x_i,y_i) + A_0 $, with $(t_i,s_i,x_i,y_i) \in A$ and 
$A_0=[0,1/4)^2 \times [0,1/2)^2$.  Using the independence of $U$ and 
$\tilde{U}$ and the scaling for the stationary pinned string with 
$L = 2^{-1/2}$, we obtain the following equality in law: 
\begin{eqnarray*}
\lefteqn{m(\tilde{V}(A_i))} \\
& = & m\left( \left\{U_{1+t}(x) - \tilde{U}_{1+s}(y) : (t,s,x,y) 
\in (t_i,s_i,x_i,y_i) + A_0 \right\} \right)   \\
& \stackrel{\mathcal{L}}{=}  & m\Big( \Big\{2^{-1/2}\left(U_{t}(x) - 
\tilde{U}_{s}(y)\right)   \\
&& \qquad : (t,s,x,y) \in (4+4 t_i, 4+ 4 s_i, 2x_i,  2y_i) + 
A) \Big\} \Big) \\
& = & \frac{1}{64} m\left( \left\{ U_{t}(x) - \tilde{U}_{s}(y) : (t,s,x,y) 
\in (4+4 t_i, 4+ 4 s_i, 2x_i,  2y_i) + A) \right\} \right) \\
& = & \frac{1}{64}  m\Big( \Big\{\Big[ U_{4+4t_i+t}(2x_i+x) - 
U_{3+4t_i}(2x_i)\Big]    \\
&& \qquad\qquad - \Big[ U_{4+4s_i+s}(2y_i+y) - U_{3+4s_i}(2y_i) \Big]
: (t,s,x,y) \in  A) \Big\} \Big)		\\
& \stackrel{\mathcal{L}}{=}  & \frac{1}{64}  m(\tilde{V}(A)). 
\end{eqnarray*}
We again obtain, using the inclusion-exclusion argument from the proof of 
theorem \ref{hitting1}, 
\[
E \left[ m \left(\tilde{V}(A_i) \cap \tilde{V}(A_j) \right) \right] =0  
\mbox{ for $i \neq j$}.  
\]
We may assume that
\[ 
A_1 = [0,1/4)^2 \times [0,1/2)^2, \quad A_2 = [1/4,1/2)^2 \times [0,1/2)^2.
\] 
Defining, for $(t,s,x,y) \in A_0$, 
\begin{eqnarray*}
V^{(1)}(t,s,x,y) & = & \left(U_{(5/4) + t}(x) -U_{5/4}(0) \right) - 
\left(\tilde{U}_{(5/4) + t}(x)
 - \tilde{U}_{5/4}(0) \right),  \\
V^{(2)}(t,s,x,y) & = & \left(U_{(5/4) - t}(x) -U_{5/4}(0) \right) - 
\left(\tilde{U}_{(5/4) - t}(x)
 -\tilde{U}_{5/4}(0) \right),
\end{eqnarray*}
we note that  
\[ 
m \left(\tilde{V}(A_i) \cap \tilde{V}(A_j) \right)= m \left(V^{(1)}(A_0) 
\cap V^{(2)}(A_0) \right).
\] 
Arguing exactly as in the proof of part (b) we may conclude that 
\[
E[m(\tilde{V}(A))] = 64 E[m(\tilde{V}(A_1))]=0, 
\]
which finishes the proof of 
(\ref{d=12key}).  \qed
%
%%%%%%%%%%%%%%%%%%%%%%%%%%%%%%%%%%%%%%%%%%%%%%%%%%%%
\section{Transience and recurrence} \label{s6}
\setcounter{equation}{0}
%
%%%%%%%%%%%%%%%%%%%%%%%%%%%%%%%%%%%%%%%%%%%%%%%%%%%%
For  the $N$-parameter Brownian sheet in $d$ dimensions, Orey and Pruitt  
\cite{op73} gave necessary and sufficient conditions on $d$ and $N$ for 
recurrence.  In this section, we will study the same question for the 
stationary pinned string $(U_t(x))$ in $\mathbf{R}^d$.  We say that a 
continuous function 
$(f_t(x): t \geq 0, x \in \mathbf{R})$ is recurrent if for any $\delta>0$ 
there exist sequences  
$(x_n), \, (t_n)$, with $\lim_{n \to \infty} t_n=\infty $, so that 
$f_{t_n}(x_n) \in B_{\delta}(0)$.  The aim of this section is to prove the 
following result.
\begin{theorem} \label{recurrence}
The stationary pinned string $(U_t(x))$ in $\mathbf{R}^d$ is almost surely 
recurrent if $d \leq 6$ and almost surely not recurrent if $d \geq 7$. 
\end{theorem}
To help in the proof of Theorem \ref{recurrence}, we first establish the 
following 0-1 law.  Define
\begin{eqnarray*}
\lefteqn{\mathcal{G}_N=  \sigma \left\{U_0(x):  |x|>N \right\}   }\\
&& \hspace{.4in} \vee \sigma \left\{W(\varphi): 
\varphi(t,x)=0 \mbox{ if $0\le t\le N$ and $|x|\le N$}\right\}.  
\end{eqnarray*}
We then set $\mathcal{G}  = \bigcap_{N=1}^{\infty} \mathcal{G}_N$.  
We can also show that $\mathcal{G}$ is trivial, using
the independence of $U_0$ and $W$, and the arguments used to prove
Kolmogorov's 0-1 law on the triviality of the Brownian tail $\sigma$-field.
\begin{lemma} \label{0-1-law}
Let $\mathcal{R}(\delta)$ be the event that there exist sequences 
$(x_n), \, (t_n)$, with $t_n \to \infty$, so that 
$U_{t_n}(x_n) \in B_{\delta}(0)$.  

Then $(\mathcal{R}(\delta): \delta>0)$ and $\mathcal{R}$ are all tail 
events in $ \mathcal{G}$.
\end{lemma}
\textbf{Proof of Lemma \ref{0-1-law}.}
For $N \geq 1$ and $t \geq N$, define 
\[
f^{(N)}_t(x)=\int_{-N}^{N}G_t(x-y)U_0(y)dy + 
\int_{0}^{N}\int_{-N}^{N}G_{t-s}(x-y)W(dy\,ds)  
\]
and set $U_t^{(N)}(x) = U_t(x) - f^{(N)}_t(x)$.  Then subtracting 
$f^{(N)}_t(x)$ from the representation for $U_t(x)$ given in 
(\ref{stationarystring}) shows that 
$(U_t^{(N)}(x): t \geq N, \, x \in \mathbf{R})$ is 
$\mathcal{G}_N$-measurable.  We claim that
\begin{equation}
\label{claim10}
\lim_{t \to \infty} \sup_{x \in \mathbf{R}} \left| f^{(N)}_t(x) \right|=0 
\end{equation}
almost surely, for each $N \geq 1$.
Assuming this claim then,  since $B_{\delta}(0)$ is an open box, we see 
that the event $\mathcal{R}(\delta)$ is unchanged, up to a null set, if we 
replace $U_t(x)$ by $U^{(N)}_t(x)$ in its definition., implying that 
$\mathcal{R}(\delta)$ is a tail event.

To prove the claim (\ref{claim10}), note that 
$ f^{(N)}_t(x) =  \int G_{t-N} (x-z) g^{(N)}(z) dz $,
where
\[ 
g^{(N)}(z) = \int_{-N}^{N}G_N (x-y) U_0(y)dy 
+\int_{0}^{N}\int_{-N}^{N}G_{N-s}(x-y) W(dy\,ds).  
\]
It is straightforward to show that $g^{(N)}(x)$ is almost surely in $L^1$. 
Then the inequality 
$ \| f_t^{(N)}\|_{\infty} \leq (4 \pi t)^{-1/2} \|g^{(N)}\|_1 $ implies the 
claim (\ref{claim10}).  \qed

\vspace{.1in}

\noindent \textbf{Proof of Theorem \ref{recurrence} in dimensions $d \leq 6$.} 
By projection, it suffices to deal with the case $d=6$.  We will use an 
inclusion-exclusion argument again, working with values of the string 
$U_t(x)$ when $t$ and $x$ are integers.  Fix $\delta \in (0,1]$, and, for 
integers $i,j$ define 
\[ 
\mathcal{R}_{i,j}= \left\{U_i(j) \in B_{\delta}(0) \right\}, \quad
\mathcal{R}(N,\delta) = \bigcup_{i=N}^{N^2} \bigcup_{0 \leq j \leq 
i^{1/2}}\mathcal{R}_{i,j}. 
\]
Our aim is to use an inclusion-exclusion argument to show that 
$ P( \mathcal{R}(N,\delta)) \geq p_0 >0$ for all 
$N$ sufficiently large.  Then, using the definition in Lemma \ref{0-1-law},  
we have 
\[
P(\mathcal{R}(\delta)) \geq P( \mathcal{R}(N,\delta) \; \mbox{infinitely often}) \geq p_0 > 0.
\]
By the zero-one law $ P(\mathcal{R} (\delta))=1$  for any $\delta>0$, which 
will complete the proof of recurrence.

The variance estimates (\ref{boundoncovariance}) on $U_t(x)$ imply that 
there exist constants $c_6,\, c_7 >0$, depending only on $\delta$, so that for
$i= 0,1,\ldots$ and $j  \in \mathbf{Z}$ with $(i,j) \neq (0,0)$
\begin{equation} \label{r-bound-1}
c_6 (i^{1/2}+|j|)^{-3}  \leq  P(\mathcal{R}_{i,j}) \leq c_7 
(i^{1/2}+|j|)^{-3}.
\end{equation}
So, for sufficiently large $N$,
\begin{eqnarray} 
\sum_{i=N}^{N^2} \sum_{0 \leq j \leq i^{1/2}} P(\mathcal{R}_{i,j}) 
& \geq &  c_6  \sum_{i=N}^{N^2} \sum_{0 \leq j \leq i^{1/2}} 
(i^{1/2}+|j|)^{-3}  \nonumber \\
& \geq & \frac{c_6}{2}  \int_{N}^{N^2}  \int_{0}^{x^{1/2}} (x^{1/2}+y)^{-3} 
dy \, dx \nonumber \\
& =  & \frac{3c_6}{16} \log(N).  \label{r-upper-2}
\end{eqnarray}
A similar calculation, using the upper bound in (\ref{r-bound-1}), shows 
that for sufficiently large $N$
\begin{equation} \label{similar-calc}
\sum_{i=N}^{N^2}\sum_{0 \leq j \leq i^{1/2}} P(\mathcal{R}_{i,j}) \leq 4 
c_7 \log (N).
\end{equation}
  From Lemma \ref{covarianceestimates1}, we have
\begin{eqnarray*}
\lefteqn{P\left( U_1(x) \in B_{\delta}(0), 
U_{1+s}(x+y) \in B_{\delta}(0) \right)      }\\
&&\hspace{.4in} \le c_3 \delta^{12} (s^{1/2} + |y|)^{-3} \quad \mbox{whenever 
$x,\, y \in [-2,2], \; s \in [0,1]$.}
\end{eqnarray*}
Using the scaling for the stationary pinned string, with the choice 
$L=t^{1/4}$, we obtain
$c_8>0$, depending only on $\delta$ so that
\begin{equation} \label{finalone}
P\left( U_t(x) \in B_{\delta}(0), U_{t+s}(x+y) \in B_{\delta}(0) \right) \leq
c_8  (t^{1/2}+|x|)^{-3} (s^{1/2} + |y|)^{-3} 
\end{equation}
whenever $t \geq 1$, $|x|,\, |y| \leq 2 t^{1/2}$ and  $s \in [0,t]$.  We 
need the bound (\ref{finalone}) for a larger set of parameters.  Since the 
stationary string is a solution to (\ref{spde}) we have
\[
U_{t+s}(x+y) = \int G_s(x+y-z) U_t(z) dz + \int^s_0 \int G_{s-r}(x+y-z) 
W(dz\,dr)
\]
so that 
\begin{eqnarray*}
\lefteqn{\mbox{Var} \bigg( U_{t+s}(x+y) - E[U_{t+s}(x+y)|\mathcal{F}_t] \bigg)  }\\
&=& \mbox{Var} \left(\int^s_0 \int G_{s-r}(x+y-z) W(dz\,dr) \right) \\
&=& C s^{1/2}.
\end{eqnarray*}
Hence 
\[
P(U_{t+s}(x+y) \in B_{\delta}(0)| \mathcal{F}_t) 
  \leq C s^{-3/2}\leq  3^3 \, C\,\cdot (s^{1/2}+|y|)^{-3}, 
\]
provided $|y| \leq 2 s^{1/2}$.  Using this we see that the bound 
(\ref{finalone}) also holds, after possibly modifying the value of  $c_8$,  
whenever $ |y| \leq 2 s^{1/2}$.

Now we can estimates the covariance term for the event 
$\mathcal{R}(N,\delta)$. 
\begin{eqnarray} 
\label{r-lower-1}
\lefteqn{ \sum_{i=N}^{N^2} \sum_{0 \leq j \leq i^{1/2}} 
\sum_{\tilde{i}=N}^{N^2} 
\sum_{0 \leq \tilde{j} \leq \tilde{i}'^{1/2}} P(\mathcal{R}_{i,j} \cap 
\mathcal{R}_{\tilde{i}, \tilde{j}})  \mathbf{1}\left( (i,j) \neq (\tilde{i}, 
\tilde{j}) \right)  } \nonumber \\
& \leq & 2 \sum_{i=0}^{N^2}\sum_{j=-\infty}^{\infty} \sum_{k=0}^{N^2} 
\sum_{\ell =-\infty}^{\infty}
P(\mathcal{R}_{i,j} \cap \mathcal{R}_{i+k,j+\ell})   \\
&& \qquad \mathbf{1} \left( (k,\ell) \neq (0,0), \; 
  |j| \leq i^{1/2}, \; |l| \leq (i+k)^{1/2} \right)   \nonumber\\
& \leq & 2 c_8  \sum_{i=0}^{N^2} \sum_{j=-\infty}^{\infty} \sum_{k=0}^{N^2}
 \sum_{\ell=-\infty}^{\infty}  (i^{1/2}+|j|)^{-3} (k^{1/2}+|\ell|)^{-3} 
\mathbf{1} \left( (k,\ell) \neq (0,0) \right).    \nonumber
\end{eqnarray}
To justify the second inequality, we note that for values of  
$k \geq i/3$ we have
$|l| \leq 2 (i+k)^{1/2} \leq 2 k^{1/2}$, and we may apply (\ref{finalone}). 
For values of $ k \leq i/3$ we have  $|l| \leq (i+k)^{1/2} \leq 2 i^{1/2}$ 
and $j \leq i^{1/2}$, and again we may  apply (\ref{finalone}).
Now we bound the double sum in (\ref{r-lower-1}), for sufficiently large $N$, 
by 
\[
C \int^{N^2}_1 \! \int_{0}^{-\infty}  (x^{1/2}+y)^{-3} dy \, dx 
\leq c_9 \left(\log(N)\right)^2
\]
where $c_9$ depends only on $\delta$.  Using (\ref{r-upper-2}), 
(\ref{similar-calc}) and (\ref{r-lower-1}) with Lemma  \ref{lower-bound}, 
we see that  $ P( \mathcal{R}(N,\delta)) \geq p_0 >0$ for sufficiently 
large $N$, completing the proof.  \qed

\vspace{.1in}

\noindent \textbf{Proof of Theorem \ref{recurrence} in dimensions $d \leq 6$.} 
It again suffices, by a projection argument, to work in dimension $d=7$. 
The strategy is to study the string along a grid of points, show that 
`recurrence on this grid' is impossible, and then to control the pieces 
between the grid points.  We define squares in the $(t,x)$ plane as follows.  
Let 
$ S_{i,j}= [i,i+1] \times [j,j+1]$ for $ i =1,2,\ldots$ and 
$j \in \mathbf{Z}$.  We will divide the squares $S_{i,j}$ into rectangles.  
To this end, let $m(i,j)$ be the unique integer such that 
\begin{equation} \label{m-def}
m(i,j)^3 \leq \left(i^{1/2}+|j| \right)^{1/4} < (m(i,j)+1)^3.  
\end{equation}
We divide each square $S_{i,j}$ into $m(i,j)^3$ rectangles, each a 
translate of $[0,m^{-2}] \times [0,m^{-1}]$, where $m=m(i,j)$.  We say these 
rectangles are of type $m$.  Let $M(m)$ be the number of  rectangles of 
type $m$, let $( R_k^{(m)}:  k=1,\dots,M(m)) $ be an enumeration of the 
rectangles of type $m$, and let $(t_k^{(m)},x_k^{(m)})$ be the point in 
$R_k^{(m)}$ with smallest $(x,t)$ coordinates.  Fix $\delta>0$.  Then, using 
the lower bound on the variance of $U_t(x)$ in (\ref{boundoncovariance}),
and the subdivision of  $S_{i,j}$ into $m(i,j)^3$ rectangles, we have
\begin{eqnarray}  
\label{borel-3}
\sum_{m=1}^{\infty}\sum_{k=1}^{M(m)} P\left(U_{t_k^{(m)}} 
\left(x_k^{(m)}\right) \in B_{2\delta}(0) \right)
& \leq & C \sum_{i=1}^{\infty} \sum_{j \in \mathbf{Z}} m(i,j)^3  
\left(i^{1/2}+|j| \right)^{-7/2}  \nonumber\\
&\leq & C \sum_{i=1}^{\infty} \sum_{j \in \mathbf{Z}}  \left( i^{1/2}+|j| 
\right)^{-13/4}	\nonumber\\
&<& \infty.  
\end{eqnarray}
The finiteness of the double sum follows by bounding it by a suitable 
integral in the usual way.  By the Borel-Cantelli lemma, the string, 
evaluated at the grid points $(t_k^{(m)}, x_k^{(m)})$, will
eventually leave the box $B_{2 \delta}(0)$.  We now interpolate between the 
grid points.  Using the boundedness of the variance of $U_t(x)$ over 
$(x,t) \in [0,1]^2$, we first apply Borel's  inequality 
for Gaussian fields (see \cite{adl90} chapter II) to find constants 
$0< c_{10},c_{11} <\infty$ so that
\[ 
P\left( \sup_{ (t,x) \in [0,1]^2 } |U_t(x)| \geq \delta \right) \leq c_1 
\exp( -c_2 \delta^2)
\]
for all $\lambda >0$. 
Now by translation invariance and then scaling we have,  for any $m \geq 1$ 
and $1 \leq k \leq M(m)$,
\begin{eqnarray} \label{estimate220}
\lefteqn{
P\left( \sup_{(t,x) \in R^{(m)}_k} |U_t(x) - U_{t^{(m)}_k} (x^{(m)}_k)|
\ge \delta\right)  }\hspace{1in}\\
& = & P\left( \sup_{(t,x) \in [0,1/m^{2}] \times [0,1/m]} |U_t(x)| 
 \geq \delta\right) \nonumber \\
& = & P\left( \sup_{(t,x) \in [0,1]^2} |U_t(x)| \geq m^{1/2} \delta\right) 
	\nonumber \\
& \leq & c_1 \exp\left( - c_2 m \delta^2\right).  \nonumber
\end{eqnarray}
We can bound for the number $M(m)$ of rectangles of type $m$ as follows.  
$M(m)$ equals $m^3$ times the number of  squares $S_{i,j}$ with $m(i,j)=m$.  
Now, (\ref{m-def}) implies that $ i \leq (m+1)^{24} $ and 
$ |j| \leq (m+1)^{12}$.  So a crude bound on $M(m)$ is given by 
$ M(m) \leq Cm^3 m^{24} m^{12} = Cm^{39}.$ Combining this with 
(\ref{estimate220}) we have
\[ 
\sum_{m=1}^{\infty} \sum_{k=1}^{M(m)}  P\left( \sup_{(t,x) \in R^{(m)}_k} 
\left|U_t(x) - U_{t^{(m)}_k} 
(x^{(m)}_k)\right| \geq \delta\right) < \infty.
\]
Combining this with (\ref{borel-3}) we may apply the Borel-Cantelli lemma 
to conclude that the probability of recurrence is zero, completing the 
proof.  \qed
%
%%%%%%%%%%%%%%%%%%%%%%%%%%%%%%%%%%%%%%%%%%%%%%%%%%%%

%
%%%%%%%%%%%%%%%%%%%%%%%%%%%%%%%%%%%%%%%%%%%%%%%%%%%%
\end{document}